\documentclass[11pt]{amsart}
\usepackage{amsmath,amssymb,amsfonts,amscd,graphicx}
\usepackage{latexsym}
\usepackage{amsmath}
\usepackage{xypic}
\usepackage[all]{xy}

\usepackage{dbnsymb}

\CompileMatrices

\numberwithin{equation}{section}

\newtheorem{cor}{Corollary}[section]
\newtheorem{prop}[cor]{Proposition}
\newtheorem{lem}[cor]{Lemma}
\newtheorem{notation}[cor]{Notation}
\newtheorem{them}[cor]{Theorem}

\theoremstyle{definition}
\newtheorem{defi}[cor]{Definition}
\newtheorem{remark}[cor]{Remark}

\newcommand{\co}{\colon\thinspace}

\begin{document}

%
%
%
%
\title[Link invariants, the chromatic polynomial and the Potts model]
{Link invariants, the chromatic polynomial\\ and the Potts model}

\author[Paul Fendley and Vyacheslav Krushkal]{Paul Fendley$^1$ and Vyacheslav Krushkal$^2$}
\address{Paul Fendley, Department of Physics, University of Virginia,
Charlottesville, VA 22904 USA;
and All Souls College and the Rudolf Peierls Centre for Theoretical Physics,
University of Oxford, 1 Keble Road, Oxford OX13NP, UK
}

\email{fendley\char 64 virginia.edu}

\address{Vyacheslav Krushkal, Department of Mathematics, University of Virginia,
Charlottesville, VA 22904-4137 USA; and Kavli Institute for Theoretical Physics, Santa Barbara, CA 93109}

\email{krushkal\char 64 virginia.edu}

\thanks{$^1$ Supported in part by NSF grants DMR-0412956 and DMR/MSPA-0704666, and by the UK EPSRC under grant EP/F008880/1}
\thanks{$^2$ Supported in part by NSF grants DMS-0729032 and PHY05-51164}

\begin{abstract}
We study the connections between link invariants, the chromatic
polynomial, geometric representations of models of statistical
mechanics, and their common underlying algebraic structure.  We
establish a relation between several algebras and their associated
combinatorial and topological quantities. In particular, we define the
{\em chromatic algebra}, whose Markov trace is the chromatic
polynomial $\chi^{}_Q$ of an associated graph, and we give
applications of this new algebraic approach to the combinatorial
properties of the chromatic polynomial.  In statistical mechanics,
this algebra occurs in the low temperature expansion of the $Q$-state
Potts model. We establish a relationship between the chromatic algebra
and the $SO(3)$ Birman-Murakami-Wenzl algebra, which is an
algebra-level analogue of the correspondence between the $SO(3)$
Kauffman polynomial and the chromatic polynomial.
\end{abstract}

\maketitle

\section{Introduction}
\label{sec:intro}

The connections between algebras, statistical mechanics, link
invariants, and topological quantum field theories have long been
exploited to great effect. The simplest and best-understood example of
such involves the Temperley-Lieb algebra, the Potts model, the Jones
polynomial, and $SU(2)$ Chern-Simons gauge theory.  The study of these
connections originated with Temperley and Lieb's work on the Potts
model in statistical mechanics \cite{TL}. They showed how to write the
transfer matrix of several two-dimensional lattice models, including
the $Q$-state Potts model, in terms of the generators of an algebra
which bears their name. Writing the transfer matrix in terms of the
generators of the TL algebra is very useful for statistical mechanics
because a number of physical quantities of the system follow purely
from the properties of this algebra, not its presentation. For
example, this rewriting yields the result that when $Q\le 4$, the
self-dual point of the Potts model is critical, whereas for $Q>4$ it
is not.

Recently such connections found an important new application in
condensed matter physics, in the study of topological states of
matter, cf \cite{F}, \cite{FFNWW}, \cite{Fnew}, \cite{LW}. In this
paper we present a number of results at the intersection of
combinatorics, quantum topology, and statistical mechanics which have
applications both in mathematics, specifically to the properties of
the chromatic polynomial of planar graphs and its relation to link
invariants, and in physics (in the study of the Potts model and of
quantum loop models).  We explain how the {\em
chromatic algebra} provides a natural setting for studying
algebraic-combinatorial properties of the chromatic polynomial.

In particular, the chromatic algebra discussed in this paper underlies
the quantum loop models discussed in \cite{FF} and further developed
in \cite{FFNWW,Fnew} (closely related models were introduced in
\cite{LW}). Quantum loop models provide lattice spin systems whose
low-energy excitations in the continuum limit are described by a
topological quantum field theory \cite{F}.  In both classical and
quantum cases, algebraic relations such as the level-rank duality
described in \cite{FK2} allow one to map seemingly different loop
models onto each other. This turns out to be quite useful in locating
critical points in loop models \cite{FJ}, a matter of great importance
for finding quantum loop models which describe topological order. In
fact, our results can be directly applied to topological quantum field
theory. We show that the chromatic algebra is associated with the
$SO(3)$ BMW algebra, implying that correlators in the $SO(3)$
topological quantum field theory can be expressed in terms of the
chromatic polynomial.

Our results can usefully be applied to both classical loop models as
well. Expressing classical loop models algebraically as described in
this paper allows one to relate loop models to other sorts of
statistical-mechanical models, such as models where the degrees of
freedom are spins or heights. In fact, the simplest application of our
results is to the same Potts model! The representation of the TL
algebra in terms of completely packed loops described in section
\ref{sec:TL} is not the only geometric representation of the Potts
model; another is generally known as the ``low-temperature''
expansion. We will explain how the chromatic algebra naturally
describes the degrees of freedom in this low-temperature expansion.

The utility of the Temperley-Lieb algebra and its generalizations
extends beyond statistical mechanics to the study of invariants of
knots and links and of topological quantum field theories.  A
graphical presentation of the TL algebra underlies the computation of
the Jones polynomial \cite{KL}. Each knot or link may be represented
as the closure of a braid, giving rise to an element of this algebra,
and evaluating the Jones polynomial for a link corresponds (up to a
normalization) to taking the Markov trace of this
element. Subsequently, Witten showed how the Jones polynomial is also
related to computations in a three-dimensional topological field
theory, Chern-Simons theory \cite{Witten} (see \cite{W} for an
exposition geared toward mathematicians).  Such computations in
Chern-Simons theory are equivalently described in two-dimensional
conformal field theory. This sequence of relations thus comes full
circle, because these conformal field theories describe the scaling
limits of two-dimensional statistical-mechanical models at their
critical points.

The purpose of this paper is twofold. One goal is to show how these
connections between link invariants, algebras, and statistical
mechanics allow us to relate seemingly different algebras and their
evaluations.  Many generalizations of the Temperley-Lieb algebra and
the Jones polynomial are now known, and Chern-Simons theories for
other representations of $SU(2)$ and for other groups are
understood. Our main focus is on the $SO(3)$ Birman-Murakami-Wenzl
algebras \cite{BW,M}, the corresponding specialization of the
Kauffman polynomial, and the $SO(3)$ TQFTs.  We will explain how
results concerning them relate to geometric models of statistical
mechanics like the Potts model.

Another purpose of this paper is to describe the chromatic algebra
(introduced for different reasons in \cite{MW}),
where the trace of an element is given by the chromatic polynomial of
an associated planar graph.  We establish a relationship between the
chromatic, $SO(3)$ BMW, and Temperley-Lieb algebras and their
traces. The chromatic algebra has a trace pairing defined in terms of
the chromatic polynomial, and we show that for $Q\geq 4$ this pairing
defines a positive-definite Hermitian product.

A nice byproduct of our analysis is that identities for the chromatic
polynomial can be extended and derived in a more transparent fashion
by utilizing the chromatic algebra. We give an algebraic proof of
Tutte's golden identity for the chromatic polynomial in a companion
publication \cite{FK2}. This striking non-linear identity plays a very
interesting role in describing quantum loop models of ``Fibonacci
anyons'', where it implies that these loop models should yield
topological quantum field theories in the continuum limit
\cite{LW,FF,FFNWW,F}. In our companion paper we also use the
Jones-Wenzl projectors in the chromatic algebra to derive linear
identities for the chromatic polynomial.

Several authors have considered similar algebraic constructions, for
example Jones \cite{J1} in the context of planar algebras, Kuperberg
\cite{Ku} in the rank 2 case, Martin and Woodcock \cite{MW} for
deformations of Schur algebras, Koo and Saleur in the setting of 
integrable lattice models (cf \cite{KooSaleur}) 
and Walker \cite{W,Walker} in the TQFT
setting.  Our approach and results are different: we derive new
relations between the chromatic and the BMW and TL algebras, and we
give applications to the structure of the chromatic polynomial of
planar graphs.

In section \ref{sec:TL}, we review the Temperley-Lieb algebra, the
Jones polynomial, and the Potts and completely packed models of
statistical mechanics.  In section \ref{sec:chromatic}, we introduce
the chromatic algebra, and show how its evaluation gives the chromatic
polynomial of the graph dual to the loop configuration. In section
\ref{TQFT section} we discuss how this algebra may be used to
construct the (doubled) $SO(3)$ topological quantum field
theory. Section \ref{sec:trivalent} establishes a presentation of the
chromatic algebra in terms of trivalent graphs.  We relate the
chromatic algebra to the $SO(3)$ BMW algebra in sections
\ref{sec:equivalence}, \ref{sec:relations}, and as a consequence show
that their evaluations are equal. A physical reason for this
equivalence has been described in depth in \cite{FR,FF}, and will be
reviewed in section \ref{sec:equivalence}.  The properties of the
trace product on the chromatic algebra are considered in section
\ref{sec:inner product}. The paper is concluded by a list of open
questions in section \ref{sec:questions}.

\section{The Temperley-Lieb algebra, the Jones polynomial and
  statistical mechanics}
\label{sec:TL}

The Temperley-Lieb (TL) algebra in degree $n$, $TL_n$, is an algebra over ${\mathbb
C}[d]$ generated by $1, e_1,\ldots, e_{n-1}$ with the relations
\begin{equation}
\label{TL relations}
e_i^2=e_i, \qquad\, e_ie_{i\pm 1} e_i=\frac{1}{d^2}\;e_i, \qquad\,
e_i e_j=e_j e_i \; \; {\rm for} \; \;
|i-j|>1.
\end{equation}
Define $TL=\cup_n TL_n$.  The indeterminate $d$ may be set to equal a
specific complex number, and when necessary, we will include the
parameter $d$ in the notation, $TL^d_n$.

\begin{figure}[ht]
\includegraphics[width=9.7cm]{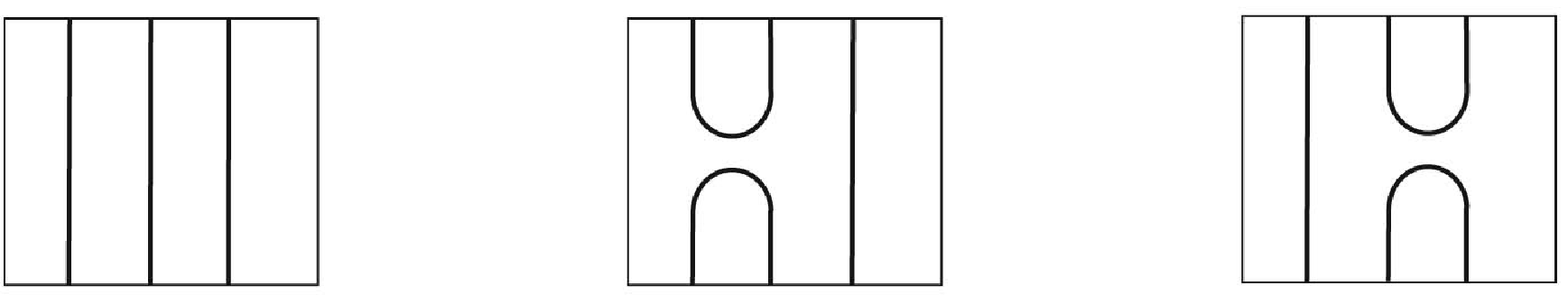}
    \put(-298,22){$1 =$}
    \put(-205,22){$e_1 =\frac{1}{d}$}
    \put(-96,22){$e_2 = \frac{1}{d}$}
\caption{Generators of $TL_3$}
\label{fig:TL}
\end{figure}

A presentation of the TL algebra of broad interest is the ``loop'' or
``$d$-isotopic'' representation, where the relations of the algebra
have a simple geometric interpretation. The loop representation of
$TL_n$ acts on a collection of $n$ strands as illustrated for $TL_3$
in figure \ref{fig:TL}.  In this setting, an element of $TL_n$ is a
formal linear combination of $1$-dimensional submanifolds in a
rectangle $R$. Each submanifold meets both the top and the bottom of
the rectangle in exactly $n$ points.  The multiplication then
corresponds to vertical stacking of rectangles.  These strands are
forbidden to cross, but we do allow adjacent strands to join, as
displayed in the figure. One can intuitively think of these as the
world lines of particles moving in one dimension; the joining of
adjacent strands corresponds to pair annihilation and recreation. The
generators $e_i$ of the TL algebra in this representation annihilate
and recreate the $i$th and $i+1$st particles.

A nice feature of this representation is that requiring that the $e_i$
satisfy the TL algebra imposes ``$d$-isotopy'', cf \cite{F}. Namely, a
circle (simple closed curve) can be removed by multiplying the
corresponding element in $TL$ by $d$.  All isotopic pictures (which
can be deformed continuously into each other without lines crossing
and while keeping the points on the boundary fixed) are considered
equivalent.

Various presentations of the TL algebra can be used to define lattice
statistical-mechanical models. When the $e_i$ are represented by
matrices, the degrees of freedom are usually referred to as spins or
heights. For example, in the $Q$-state Potts model, the degrees of
freedom are ``spin'' variables $\sigma_i$ taking integer values
$1\dots Q$ at each site of some lattice. The $e_i$ here are
represented by tensor products of $Q\times Q$ matrices, with
$d=\sqrt{Q}$. The transfer matrix of the $Q$-state Potts model, with
$n$ the number of sites on a zig-zag line, can be written entirely in terms of
these generators.  At the isotropic self-dual point for
$n$ even, this transfer matrix for the square lattice is
\begin{equation}
T=\left[\prod_{j=1}^{n/2} (1+e_{2j})\right] \left[\prod_{j=1}^{n/2} (1+e_{2j-1})\right]
\label{transfer}
\end{equation}
The partition function for an $n\times m$ system with periodic
boundary conditions in the $m$ direction is then $Z=\hbox{tr}\, T^m$.

A closely related way of defining a statistical mechanical model is
via a pictorial presentation. For this presentation of the TL
algebra, this results in the {\em completely packed loop model}.  This
model is defined on any graph with four edges per vertex. Each edge of
the graph is covered by a loop, and at each vertex the loops avoid
each other in the two possible ways. By using the pictures in figure
\ref{fig:TL}, each configuration then corresponds to a single word in
the TL algebra (i.e.\ some product of the $e_i$). The transfer matrix
on the square lattice at the self-dual point remains (\ref{transfer});
each factor $(1+e_i)$ describes the two choices for the loops'
behavior at a single vertex. Expanding the product in $T^m$ into
individual words corresponds to writing the partition function as a
sum over loop configurations. The Boltzmann weight of each
configuration in the completely packed loop model is the
``evaluation'' or the ``Markov trace'' of the corresponding element.
This is a linear map of the algebra to the complex numbers, namely,
the trace $tr_d\co TL^d_n\longrightarrow {\mathbb C}$ is defined on
the additive generators (rectangular pictures) by connecting the top
and bottom endpoints by disjoint arcs in the complement of the rectangle $R$ in the
plane. The result is a disjoint collection of circles in the plane,
which are then evaluated by taking $d^{\# circles}$.
This completely packed loop model is still often referred to
as the Potts model in its ``Fortuin-Kasteleyn'' or ``cluster'' representation
\cite{FK}. Note however that although the Potts model is originally
defined with $Q=d^2$ an integer, in this representation of the TL
algebra this constraint is no longer required.

Utilizing this pictorial representation and the Markov trace allows
one to relate the Jones polynomial of knot theory to the TL algebra
\cite{KL}. Namely, the Jones polynomial can then be computed by
projecting a knot or a link onto the plane, where in non-trivial cases
the projection will include overcrossings and undercrossings.  These
are described by another geometric realization of $TL_n$.  Here the
elements are framed tangles in $D^2\times [0,1]$ which meet the top
and the bottom of the cylinder in $n$ specified points, modulo the isotopy and the
skein relations
\begin{equation} \label{Jones relation} \slashoverback \; = \;
\; A\; \smoothing \; + \; A^{-1}\; \hsmoothing \qquad{\rm and}
\qquad {\mathcal L}\; \cup \bigcirc \; =\; (-A^2-A^{-2})\; {\mathcal L}\ .
\end{equation}
We call this the skein-theoretic version, $TL_n^{\rm skein}$, while
the planar one discussed above is $TL_n^{\rm planar}$.
The planar approach is more suitable for applications
to $2 D$ lattice models, while skein theory provides a more well-known
route to the construction of topological quantum field theories.
However, it is easy to see that if we set
$d=-A^2-A^{-2}$, one has the isomorphisms $TL_n\cong TL_n^{\rm skein}\cong TL_n^{\rm planar}$,
and we use the superscript only to indicate a specific geometric
context.

Up to overall factors of $A$, the Jones polynomial for a given
collection of links is the Markov trace of the corresponding
element of the algebra.  Precisely, the trace $tr_{\langle\rangle}\co
TL_n^{\rm skein}\longrightarrow {\mathbb C}$ is defined on
the generators (framed tangles) by connecting the top and bottom
endpoints by standard arcs in the complement of $D^2\times[0,1]$ in
$3$-space, sweeping from top to bottom,
and computing the Kauffman bracket \cite{KL}. (The Kauffman bracket $\langle
L\rangle$ of a link $L$ in $S^3$ is defined by the skein relations
(\ref{Jones relation}).)  It follows from these definitions that the two
traces on $TL_n$ are equal up to the change of basis, in other words
the diagram
\begin{equation} \label{TL traces} \xymatrix{ TL_n^{\rm planar}  \ar[d]^{tr_d} \ar[r]^{\cong} & TL_n^{\rm skein}  \ar[d]^{tr_{\langle\rangle}}\\
{\mathbb C}\ar[r]^=  &  {\mathbb C} }
\end{equation}
commutes. The isomorphism above is given by viewing the generators of $TL_n^{\rm planar}$
as elements in
$TL_n^{\rm skein}$ by including the rectangle $R$ as a vertical slice
of the cylinder $D^2\times[0,1]$. The inverse map $TL_n^{\rm
skein}\longrightarrow TL_n^{\rm planar}$ is defined by resolving any
tangle, using the skein relation (\ref{Jones relation}), into a linear
combination of embedded planar pictures.

The relation of the foregoing to topological field theory is well
known: the Jones polynomial of a collection of links correspond to a
correlation function of Wilson loops in $SU(2)$ Chern-Simons gauge
theory \cite{Witten,W}. These Wilson loops transform in the spin-1/2
representation of the $SU(2)$ algebra.  Since the spin $1/2$
representation of $SU(2)$ is the simplest non-trivial representation
of the simplest non-abelian Lie algebra, it is natural to expect that
the previous results have a myriad of generalizations.

\section{The chromatic algebra} \label{sec:chromatic}

The purpose of this section is to define the ``chromatic algebra''. It
is studied in further detail in section \ref{sec:trivalent}, and in
sections \ref{sec:equivalence}, \ref{sec:relations} the chromatic
algebra is related to the $SO(3)$ BMW algebra.

The {\em chromatic polynomial} ${\chi}_{\Gamma}(Q)$ of a graph
$\Gamma$, for $Q\in {\mathbb Z}_+$, is the number of colorings of the
vertices of $\Gamma$ with the colors $1,\ldots, Q$ where no two
adjacent vertices have the same color. To study ${\chi}_{\Gamma}(Q)$
for non-integer values of $Q$, it is often convenient to utilize the
{\em contraction-deletion relation} (cf \cite{B}). Given any edge $e$
of $\Gamma$ which is not a loop,
\begin{equation}
\label{chromatic poly1}
\chi_{\Gamma}^{}(Q)={\chi}_{{\Gamma}\backslash e}(Q)-{\chi}_{{\Gamma}/e}(Q)
\end{equation}
where ${\Gamma}\backslash e$ is the graph obtained from $\Gamma$ by
deleting $e$, and ${\Gamma}/e$ is obtained from $\Gamma$ by
contracting $e$. (If $\Gamma$ contains a loop then
${\chi}_{\Gamma}\equiv 0$.) A useful consequence of the
contraction-deletion relation is that
\begin{equation} \label{chromatic poly2}
{\chi}_{\Gamma}(Q)=\sum_{S\subset\{{\rm edges \;\, of}\;\,  {\Gamma}\}} (-1)^{|S|}\,
Q^{k(S)}
\end{equation}
where $k(S)$ is the number of connected components of the graph which
has the same vertices as ${\Gamma}$ and whose edge set is given by
$S$. Either of these two equations, together with the value on the
graph consisting of a single vertex and no edges: ${\chi}(\cdot)=Q$
determines the chromatic polynomial, and may be used to define it for
any (not necessarily integer) value of $Q$.


In section \ref{sec:TL} we described how the completely packed loop model
is related to (and can be defined using) the Temperley-Lieb
algebra. To motivate what follows, it is useful to describe the
analogous lattice statistical-mechanical model here. Remarkably, this
model is also closely related to the Potts model, just like the
completely packed loops. Instead of the FK/cluster expansion,
the geometric degrees of freedom of interest here arise in the
low-temperature expansion.

Each configuration of the Potts model is given by specifying the value
$\sigma_i=1\dots Q$ of a spin at each vertex of some graph $L$ (which
in physics applications is typically a lattice, but need not be). The
Boltzmann weight of each configuration is then $e^{-\beta {\mathcal
E}}$, where $\beta$ is inverse temperature, and the energy is
\begin{equation}
{\mathcal E} = - J\sum_{<ij>} \delta_{\sigma_i\sigma_j}
\label{EPotts}
\end{equation}
for nearest-neighbor sites labeled by $i$ and $j$. $J$ is a coupling,
so that when $J>0$ the model is ferromagnetic, and when $J<0$ it is
antiferromagnetic. The
partition function is then defined as
\begin{equation}
Z=\sum_{\{\sigma_i = 1,\dots,Q\}} e^{-\beta {\mathcal E}}
\label{ZPotts}
\end{equation}
where the sum is over all configurations.

To understand the low-temperature expansion, it is useful to first
describe the zero-temperature antiferromagnetic limit, where
$\beta\to\infty$ and $J<0$. The only configurations which contribute
to the sum for $Z$ in (\ref{ZPotts}) in this limit are those in which
adjacent spins have different values. $Z$ then simply counts the
number of such configurations, because each has the same weight $1$.
Thus when $J<0$,
$$\lim_{\beta\to\infty} Z = \chi_L(Q) \ .$$ Thus the chromatic
polynomial arises very naturally in statistical mechanics.

\begin{figure}[ht]
\vspace{.2cm}
\includegraphics[width=2.5cm]{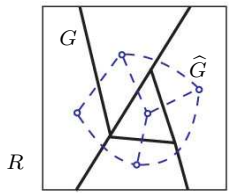}
{\scriptsize
    \put(-85,10){$R$}
    \put(-65,57){$G$}
    \put(-16,45){$\widehat G$}}
\vspace{.2cm} \caption{An element $G$ of ${\mathcal G}_2$ and its dual graph
$\widehat G$ (drawn dashed).}
\label{fig:gghat}
\end{figure}

The low-temperature expansion is an expansion of $Z$ in powers of
$e^{\beta J}$. It is useful intuitively to describe this in terms of
{\em domain walls} on the dual graph $\widehat{L}$.  Given
$G\in{\mathcal G}$, the vertices of its {\em dual graph} $\widehat G$
correspond to the complementary regions $R\smallsetminus G$, and two
vertices are joined by an edge in $\widehat G$ if and only if the
corresponding regions share an edge, as shown in figure
\ref{fig:gghat}.

With each
configuration of spins, one associates a subgraph $N$ of $\widehat{L}$
by the following rule: when the spins on two sites differ, then the
edge separating them belongs to $N$. If the spins are
the same, the corresponding edge is not part of $N$. The graph $N$ is
a domain-wall configuration, separating domains of like spins from
each other.  By construction, a domain-wall configuration $N$ consists
of a graph with no ends (no 1-valent vertices) except possibly at the
outside of $L$. For this reason we call such graphs ``nets''.  In
the zero-temperature limit, all spin configurations contributing to
the sum in $Z$ are associated to a single net $N=\widehat{L}$.

The idea of the low-temperature expansion is to first find the
domain-wall configuration $N$ associated to a given spin
configuration. Typically, many spin configurations
are associated to the same $N$. By construction, the number of these is
precisely the chromatic polynomial $\chi_{\widehat{N}}(Q)$. Each of
these has the same Boltzmann weight $e^{\beta J
(E(\widehat{L})-E(N))}$, where $E(G)$ is the number of edges in the
graph $G$. $E(N)$ can be thought of as the ``length'' of the domain
walls.  The partition function is then
\begin{equation}
Z=e^{\beta J E(\widehat{L})}
\sum_{N} e^{-\beta J E(N)} \chi_{\hat{N}}(Q)
\label{Zchrome}
\end{equation}
where the sum is over all subgraphs $N$ of $\widehat{L}$.  (Here we
can ignore the restriction that $N$ have no 1-valent vertices because
$\chi_{\widehat{G}}(Q)=0$ for any graph $G$ with 1-valent vertices.)
This shows that the $Q$-state Potts model has a very natural
description as a sum over geometric objects, nets. Note also that this
allows the model defined for any $Q$.

In section \ref{sec:TL}, we explained how the partition function of the
completely packed loop model is defined as a sum over geometric
objects, with each configuration associated with an element of the
Temperley-Lieb algebra. It is thus natural to define an algebra whose
elements correspond to nets, and whose Markov trace gives the
chromatic polynomial.  We therefore define the {\em chromatic algebra}
in the same fashion as the TL defined in section \ref{sec:TL}.
Consider the set ${\mathcal G}_n$ of the isotopy classes of planar
graphs $G$ embedded in the rectangle $R$ with $n$ endpoints at the top
and $n$ endpoints at the bottom of the rectangle.  The intersection of
$G$ with the boundary of $R$ consists precisely of these $2n$ points.
In the Potts model, a graph $G$ is comprised of the domain walls
separating regions of like spins from each other.  It is convenient to
divide the set of edges of $G$ into {\em outer} edges, i.e. those
edges that have an endpoint on the boundary of $R$, and {\em inner}
edges, whose vertices are in the interior of $R$. (Note that the
graphs $G$ are not necessarily connected.) It is convenient to allow
$G$ to have connected components which are simple closed curves (which
are not strictly speaking ``graphs'' since they do not contain a
vertex.)

\begin{notation} While discussing the chromatic algebra, we will interchangeably use two variables, $Q$ and $q$.
Set $Q=q+2+q^{-1}=(q^{1/2}+q^{-1/2})^2$.
\end{notation}

The defining contraction-deletion rule (\ref{chromatic poly1}) may be
viewed as a linear relation between the graphs $G, G/e$ and
$G\backslash e$, so in this context it is natural to
consider the vector space defined by graphs, rather than just the set
of graphs.  Thus let ${\mathcal F}_n$ denote the free
algebra over ${\mathbb C}[Q]$ with free additive generators
given by the elements of ${\mathcal G}_n$.  As usual, the multiplication
is given by vertical stacking, and we set ${\mathcal F}=\cup_n {\mathcal
F}_n$.

\begin{figure}[ht]
\includegraphics[height=2cm]{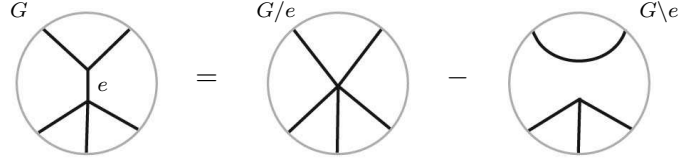}
    \put(-175,28){$=$}
    \put(-80,28){$-$}
{\scriptsize
    \put(-212,26){$e$}
    \put(-245,54){$G$}
    \put(-152,54){$G/e$}
    \put(-7,54){$G\backslash e$}}
\caption{Relation (1) in the chromatic algebra}
\label{fig:chrome1}
\end{figure}

\begin{figure}[ht]
\includegraphics[height=2cm]{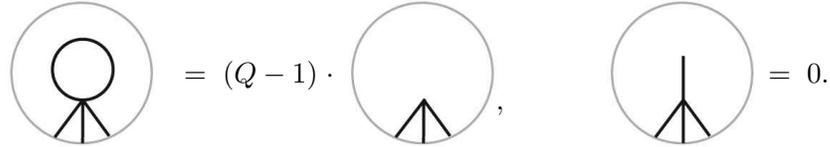}
    \put(3,25){$=\; 0.$}
    \put(-218,25){$=\; (Q-1)\;\cdot$}
    \put(-100,15){$,$}
\caption{Relations (2), (3) in the chromatic algebra}
\label{fig:chrome2}
\end{figure}

The local relations among the elements of ${\mathcal G}_n$, analogous to
contraction-deletion rule for the chromatic polynomial, are given in figures \ref{fig:chrome1},
\ref{fig:chrome2}. Note that these relations only apply to
{\em inner} edges which do not connect to the top and the bottom of
the rectangle. They are

(1) If $e$ is an inner edge of a graph $G$ which is not a loop, then
$G=G/e-G\backslash e$, figure \ref{fig:chrome1}.

(2) If $G$ contains an inner edge $e$ which is a loop, then $G=(Q-1)\;  G\backslash e$,
figure \ref{fig:chrome2}. (In particular, this
relation applies if $e$ is a simple closed curve not connected to the rest of the graph.)

(3) If $G$ contains a $1$-valent vertex (in the interior of the rectangle)
then $G=0$, figure \ref{fig:chrome2}.

\medskip

\begin{defi} \label{chromatic definition} The {\em chromatic algebra} in degree $n$, ${\mathcal C}_n$, is an algebra over ${\mathbb C}[Q]$
which is defined as the quotient of the free algebra ${\mathcal F}_n$ by the ideal
$I_n$ generated by the relations (1), (2), (3). Set ${\mathcal C}=\cup_n {\mathcal C}_n$.
\end{defi}

\medskip

The ideal $I_n$ in the definition above is generated by linear combinations of graphs in ${\mathcal F}_n$
which are identical outside a disk embedded in the rectangle, and which differ according to one of the relations $(1)-(3)$ in the disk.

\begin{remark} Recall that $e$ is a {\em bridge} if
it is an internal edge which, if removed, disconnects $G$ (considering
all points on the boundary of $R$ to be connected).
The relation (3) above can be replaced by

(3$'$): If $G$ has a bridge $e$ then $G=0$.

 Note that this means the dual graph $\widehat G$ contains a loop.
\end{remark}

\begin{figure}[ht]
\vspace{.2cm}
\includegraphics[height=1.9cm]{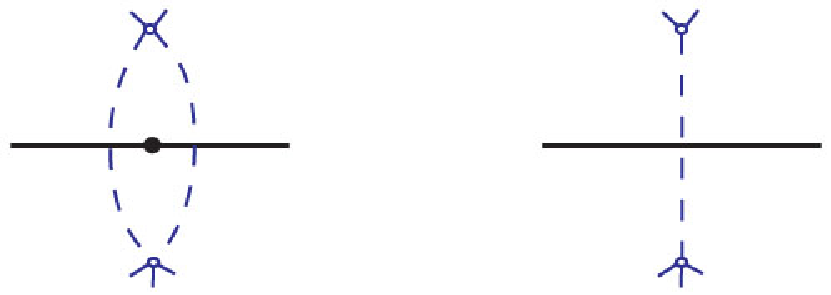}
    \put(-95,25){$=$}
\vspace{.2cm} \caption{}
\label{fig:chrome3}
\end{figure}

\begin{figure}[ht]
\vspace{.2cm}
\includegraphics[width=13cm]{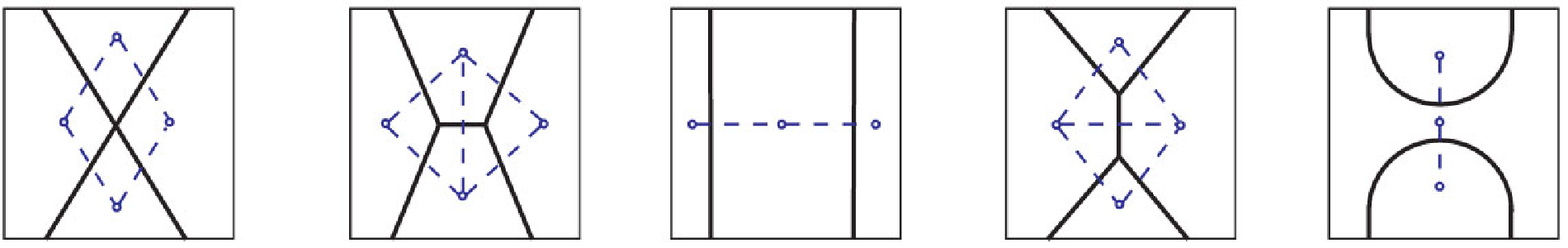}
    \put(-307,27){$=$}
    \put(-227,27){$+$}
    \put(-150,27){$=$}
    \put(-72,27){$+$}
    \vspace{.2cm}
\caption{}
\label{fig:chromeexample}
\end{figure}

We will now collect some elementary consequences of the relations which hold in ${\mathcal C}$: (1) and (3) imply that
a $2$-valent vertex may be deleted, and the two adjacent edges merged, figure \ref{fig:chrome3} (the dual graph, discussed in more detail below, is drawn dashed.) (1) and (3) also imply that if
a graph $G$ contains an isolated vertex $v$, then $G=G\backslash v$. (Note that deleting an isolated vertex does not
change the dual graph.)
Figure \ref{fig:chromeexample} gives more examples of relations which hold in ${\mathcal C}_2$. It is important to note that
the relations (1) -- (3) are
consistent with the relations for the chromatic polynomial of
the dual graph, see the following proposition.

\begin{prop} \label{well defined} \sl The chromatic polynomial gives rise to a well-defined linear map
${\chi}\co {\mathcal C}_n\longrightarrow {\mathbb C}[Q]$. This map is defined on the additive
generators $G$ as the chromatic polynomial of the dual graph, ${\chi}_{\widehat G}$, and it is
extended to ${\mathcal C}_n$ by linearity.
\end{prop}

\begin{figure}[ht]
\vspace{.2cm}
\includegraphics[width=8cm]{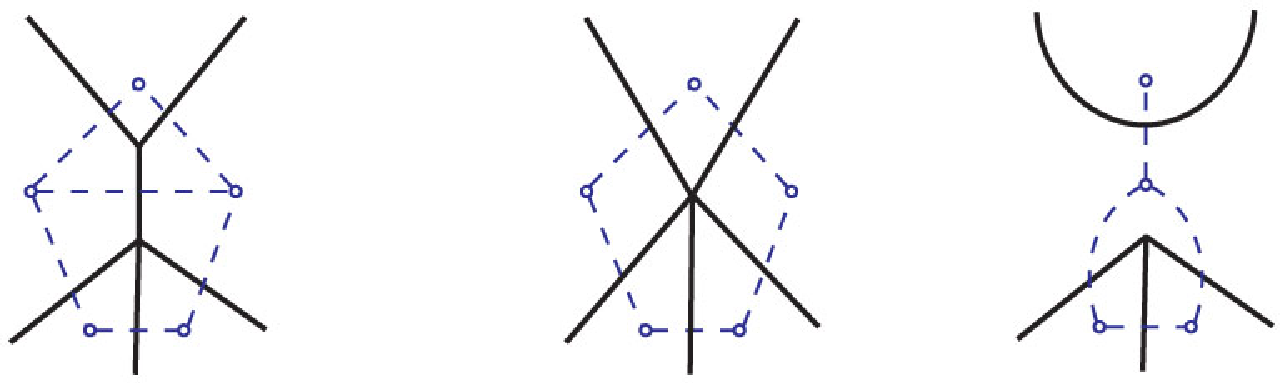}
    \put(-160,30){$=$}
    \put(-65,30){$-$}
{\scriptsize
    \put(-202,27){$e$}
    \put(-215,28){$\hat e$}
    \put(-235,60){$G$}
    \put(-144,60){$G/e$}
    \put(-23,38){$\hat e_1$}
    \put(0,60){$G\backslash e$}}
\vspace{.2cm} \caption{Relation (1) and the dual graphs}
\label{fig:relation 1 and dual}
\end{figure}

To prove this proposition, one needs to check that the relations
(1)-(3) hold when one considers the chromatic polynomial of the dual
graph. Specifically, in case (1) consider the edge $\hat e$ of
$\widehat G$, dual to $e$, figure \ref{fig:relation 1 and dual}.  Then
$\widehat{G/e}=\widehat G\backslash \hat e$, and $\widehat{G\backslash
  e}=\widehat G/ \hat e$.  (In the case when one of the vertices of
$e$ is trivalent, as in figure \ref{fig:relation 1 and dual},
$\widehat G/ \hat e$ differs from $\widehat{G\backslash e}$ by the
addition of an edge parallel to $\hat e_1$, figure \ref{fig:relation 1
  and dual}. Because of the equality in figure \ref{fig:chrome3},
${\chi}_{\widehat G/ \hat e}={\chi}_{\widehat{G\backslash e}}$ still
holds in this case.)  Therefore the relation (1) translates to
${\chi}_{\widehat G}={\chi}_{\widehat G\backslash \hat
  e}-{\chi}_{\widehat G/ \hat e}$, the defining contraction-deletion
relation for the chromatic polynomial.

To check (2), note that if $G$ contains a loop which is trivial in the plane (i.e. the region $D$
bounded by it is disjoint from the graph), the dual graph $\widehat G$ contains a
$1$-valent vertex (corresponding to the region $D$.) The effect of
deleting such a vertex and the adjacent edge on the chromatic polynomial ${\chi}_{\widehat G}(Q)$
is multiplication by $(Q-1)$. The general case when the region $D$ bounded by the
loop contains other vertices or edges of $G$ follows by applying (1) and (3) to the part
of the graph inside $D$, and then inductively applying the case of (2) considered above
to the trivial inner-most loops of $G$.

The relation (3) holds since in this case
the dual graph has a loop, therefore its chromatic polynomial vanishes.
\qed

\begin{defi} \label{chromatic trace def}
The trace, $tr_{\chi}\co {\mathcal C}\longrightarrow{\mathbb C}$ is defined on the additive
generators (graphs $G$ in the rectangle $R$) by connecting the top and bottom
endpoints of $G$ by disjoint arcs in complement of $R$ the plane (denote the result by $\overline G$)
and evaluating the chromatic polynomial of the dual graph:
\begin{equation} \label{chromatic trace}
tr_{\chi}(G)\; \, =\; \, Q^{-1}\cdot {\chi}_{\widehat{\overline G}}(Q).
\end{equation}
\end{defi}

\begin{figure}[ht]
\vspace{.2cm}
\includegraphics[height=6cm]{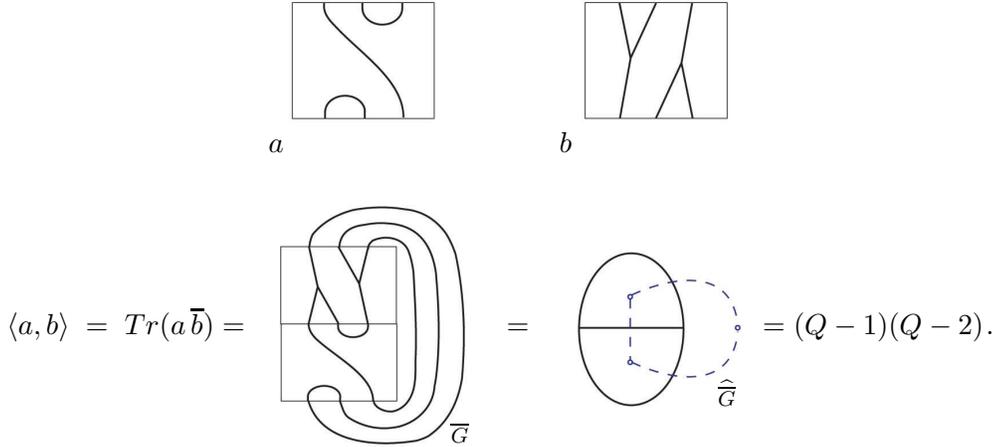}
    \put(-180,113){$a$}
    \put(-70,113){$b$}
    \put(-279,44){$\langle a,b\rangle \; =\; Tr(a\, \overline b)=$}
    \put(-90,44){$=$}
{\scriptsize
    \put(-111,3){$\overline G$}
    \put(-10,17){$\widehat{\overline G}$}}
    \put(7,44){$=(Q-1)(Q-2)$.}
 \vspace{.2cm}
\caption{The inner product in ${\mathcal C}_3$: $Tr(a\,\overline b)$ equals the evaluation
of $\overline G$, or equivalently $Q^{-1}$ times the chromatic polynomial of the dual graph
$\widehat{\overline G}$. In this
example, $\overline G$ is the theta graph, and $$\langle a,b\rangle\, =\,Q^{-1}\cdot(Q(Q-1)(Q-2))=(Q-1)(Q-2).$$}
\label{fig:innerproduct}
\end{figure}

Proposition \ref{well defined} shows that the trace is well defined.
The factor
$Q^{-1}$  is
a convenient normalization which makes the relation with the BMW algebra easier to state (see sections \ref{sec:equivalence},
\ref{sec:relations}); with this normalization, $tr_{\chi}(\cdot)=1$.
The {\em Hermitian product} on ${\mathcal C}_n$ is defined analogously to the Temperley-Lieb case:
$\langle a,b\rangle = \hbox{tr}(a\, \bar b)$, where the
involution $\bar b$ is defined by conjugating the complex
coefficients, and on an additive generator $b$ (a graph in $R$) it is
defined as the reflection in a horizontal line, see figure \ref{fig:innerproduct}.

\begin{remark} \label{Tutte specialization}
The chromatic polynomial ${\chi}_G$ and the Potts-model partition
(\ref{Zchrome}) function are specializations of the more general Tutte
polynomial $T_G(X,Y)$, cf \cite{B}. The Tutte polynomial satisfies the
well-known duality $T_G(X,Y)=T_{\widehat G}(Y,X)$, where $G$ is a
planar graph and $\widehat G$ is its dual. Therefore, rather then
defining the trace of the chromatic algebra as the chromatic
polynomial of the dual graph $\widehat G$, one could define the trace
as the corresponding relevant specialization of the Tutte polynomial
(sometimes known as the {\em flow polynomial}) of the graph $G$
itself.
\end{remark}

\begin{remark}
One could generalize the definition of ${\mathcal C}$ and consider the
``Tutte algebra'' whose trace is given by the Tutte polynomial. We
restrict our discussion to the special case of the chromatic
polynomial since it corresponds to our main object of interest, the
$SO(3)$ BMW algebra, see sections \ref{sec:equivalence},
\ref{sec:relations}. Another reason for this is that the chromatic
algebra gives rise to the $(2+1)$-dimensional $SO(3)$ TQFT (see
section \ref{TQFT section}), and this seems to be the only
specialization of the Tutte polynomial which yields a finite
dimensional unitary TQFT.
\end{remark}

\subsection{From algebras to TQFTs} \label{TQFT section}

There is a well-understood route for representing $(2+1)-$
dimensional topological quantum field theories
in terms of ``pictures'' on surfaces modulo local relations, see \cite{W}, \cite{FNWW}.
This is described in depth in the case of (doubled) $SU(2)$ theories in \cite{F}.
Here we briefly sketch the analogous construction of the doubled (Turaev-Viro \cite{TV})
$SO(3)$ TQFTs (developed in \cite{RT}, \cite{Turaev}.) In fact, the problem of
finding a description of these TQFTs in terms of ($2-$dimensional) pictures
and relations was a starting point for our introduction of the chromatic
algebra in this paper. The description of TQFTs in such terms is important
for applications in physics, specifically to lattice models exhibiting topological
order \cite{F}, \cite{LW}.

Given a compact surface ${\Sigma}$, consider the (infinite-dimensional) complex vector space $V_1$ consisting
of formal linear combinations of the isotopy classes of graphs embedded in $\Sigma$. If $\Sigma$ has
a non-empty boundary, one fixes a boundary condition -- a finite number of points in the boundary
$\partial \Sigma$, and the graphs in $\Sigma$ should meet the boundary in these specified points.
Consider the quotient $V_2^Q$ of $V_1$ by the local relations
$(1)$-$(3)$ defining the chromatic algebra as in \ref{chromatic definition}. Equivalently, one may start with the space of {\em trivalent}
graphs, modulo the relations in figure \ref{fig:trivalent}, see theorem \ref{trivalent lemma} in
the next section. The chromatic algebra ${\mathcal C}^Q$ (more precisely, its generalization, the chromatic
{\em category}, see section 4 in \cite{FK2}) is a local version of $V_2^Q$, i.e. it corresponds to ${\Sigma}=$disk.

The vector space $V_2^Q$ is still infinite-dimensional. However, at
the special values of $Q$ known as Beraha numbers,
$Q=B_n=2+2\cos(2{\pi}/n)$, there are additional local linear relations
such that the quotient, $V^{\Sigma}_n$, is finite dimensional. These
additional local relations are generators of the trace radical of the
chromatic algebra, i.e. elements $a$ of ${\mathcal C}^{B_n}$ such that
$\langle a, b\rangle =Tr(a\bar b)=0$ for all $b\in{\mathcal
C}^{B_n}$. The trace pairing descends to a positive-definite Hermitian
product on the quotient of ${\mathcal C}^{B_n}$ by the trace radical
(see corollary \ref{cor:inner product}), which gives rise to a
positive-definite Hermitian product on $V^{\Sigma}_n$, see
\cite{W}. (The trace radical is analyzed in \cite{FK2}, where we show
that it contains the pull-back of the Jones Wenzl-projector from the
Temperley-Lieb algebra, and give applications of this fact to the
structure of the chromatic polynomial of planar graphs.)

These finite-dimensional vector spaces, $V^{\Sigma}_n$, are the doubled $SO(3)$ TQFTs, and
the unitary structure on these TQFTs is induced by the trace pairing on the chromatic algebra as
indicated above. The structure of the simplest non-trivial TQFT which arises from
this construction, corresponding to $Q=B_5={\phi}+1$, and which is known as the ``doubled Fibonacci
theory'', is considered in \cite{FFNWW}. In this paper
we analyze some of the algebraic structure underlying these TQFTs for other levels $n$.

\section{A trivalent presentation of the chromatic algebra} \label{sec:trivalent}

In this section we define an algebra using {\em trivalent} graphs
modulo certain simple relations consistent with the
contraction-deletion rule, and prove that it is isomorphic to the
chromatic algebra. This result should be compared with section
\ref{sec:equivalence} which shows that, in contrast, a presentation of
this algebra in terms of four-valent graphs is rather involved. As
part of the proof, in this section we find an additive basis of the
chromatic algebra ${\mathcal C}_n$ in terms of planar partitions, and
we establish an algebra analogue of the ``state sum'' formula
(\ref{chromatic poly2}) for the chromatic polynomial.  Both the
chromatic algebra and its trivalent presentation established here are
used in our companion paper \cite{FK2} to prove and generalize Tutte's
identities for the chromatic polynomial.

\begin{figure}[ht]
\includegraphics[height=1.85cm]{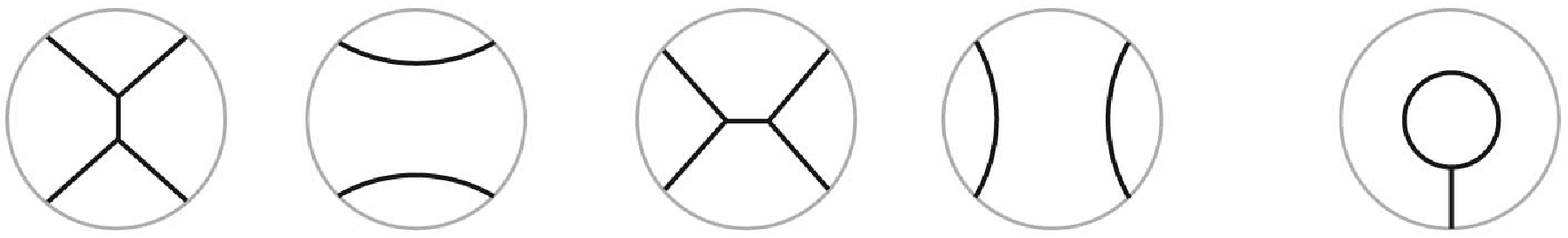}
    \put(-296,25){$+$}
    \put(-226,25){$=$}
    \put(-155,25){$+$}
    \put(4,25){$=\;0.$}
    \put(-87,8){$,$}
\caption{Trivalent presentation of the chromatic algebra.}
\label{fig:trivalent}
\end{figure}

\begin{defi} \label{trivalent definition} Analogously to definition \ref{chromatic definition},
consider the free algebra ${\mathcal {FT}}_n$ over ${\mathbb C}[Q]$
whose elements are formal linear combinations of the isotopy classes of trivalent graphs in a
rectangle $R$. The intersection of each such graph with the boundary
of $R$ consists of precisely $2n$ points: $n$ points at the top and
the bottom each.
Let ${\mathcal T}_n$ denote the quotient of ${\mathcal {FT}}_n$ by the ideal
generated by the local relations shown in
figure \ref{fig:trivalent}.
\end{defi}

\begin{remark}
Note that the vertices of the
graphs in the definition above in the interior of $R$ are trivalent, in particular they do not have ends
($1$-valent vertices) other than those on the boundary of $R$.
It is convenient
to allow $2$-valent vertices as well, so there may be loops disjoint from the rest
of the graph.

The relations in ${\mathcal T}$, shown in figure \ref{fig:trivalent},
have a rather natural interpretation in the context of TQFTs
discussed in section \ref{TQFT section}. The second relation says
that the vector space associated to the disk with a single boundary
label is trivial. The first relation implies that the vector space
associated to the disk with four boundary labels is $3$-dimensional.
\end{remark}

\begin{them} \label{trivalent lemma} \sl The map ${\Phi}\co{\mathcal T}_n\longrightarrow {\mathcal C}_n$,
induced by the inclusion of the trivalent graphs in the set of all graphs, is an algebra isomorphism.
\end{them}

First observe that $\Phi$ is well-defined. Indeed, the first
relation in figure \ref{fig:trivalent} is a consequence of the contraction-deletion rule (1) in the definition
of ${\mathcal C}_n$, see figure \ref{fig:chromeexample}. The relation on the right in figure
\ref{fig:trivalent} is a consequence of relations (2), (3) defining ${\mathcal C}_n$.
Moreover, $\Phi$ is a surjective map: using the contraction-deletion rule (1), any graph $G\in {\mathcal G}_n$
may be expressed as a linear combination of trivalent graphs. This establishes
\begin{equation} \label{dim inequality}
dim_{\mathbb C}({\mathcal T}_n)\geq dim_{\mathbb C}({\mathcal C}_n).
\end{equation}

We start the proof of the converse inequality by presenting an algebra analogue
of the expansion (\ref{chromatic poly2}) of the chromatic polynomial. This expansion is then used
for describing a linearly independent set of additive generators of ${\mathcal C}_n$.

\begin{figure}[ht]
\vspace{.2cm}
\includegraphics[height=2cm]{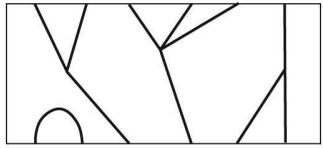}
 \vspace{.2cm}
\caption{An example of a graph in $B_6$.}
\label{fig:basis}
\end{figure}

Following the terminology introduced in section \ref{sec:chromatic},
given a graph $G$ in the rectangle $R$, $G\in {\mathcal G}_n$, we
consider its set of {\em inner} edges, the edges of $G$ whose
endpoints are not on the boundary of $R$. Consider the set $B_n$ of
all graphs in ${\mathcal G}_n$ without inner edges, as in figure
\ref{fig:basis}. (Note that the elements of $B_n$ are in $1-1$
correspondence with the {\em planar partitions} of the set of $2n$
boundary vertices such that each block of the partition contains at
least two vertices.)
\begin{lem} \label{basis lemma} \sl The elements of $B_n$ form an additive basis of the chromatic algebra ${\mathcal C}_n$.
\end{lem}

{\em Proof.} Consider the vector space $V_n$ over ${\mathbb C}$ spanned by $B_n$.
Define linear maps ${\phi}\co V_n\longrightarrow {\mathcal C}_n$
and ${\psi}\co {\mathcal C}_n\longrightarrow V_n$.
Here ${\phi}$ is induced by the inclusion $B_n=\{$isotopy classes of graphs without inner edges in the rectangle $R\} \subset
\{$isotopy classes of all graphs in $R\}={\mathcal G}_n$,
while $\psi$ is defined using an expansion similar to
the expansion (\ref{chromatic poly2}) for the chromatic polynomial: for $G\in {\mathcal G}_n$, set

\begin{equation} \label{chromatic expansion}
{\psi}(G)=\sum_{S\, \subset\, \{{\rm inner \;\,edges \;\, of}\;\,  G\}}\; (-1)^{E(G)-|S|}\,\;
Q^{n(S)}\; b_S\ ,
\end{equation}

where $E(G)$ is the number of edges of the graph $G$ and
$n(S)$ is the {\em nullity} of the graph $G_S$ which is obtained by keeping all vertices and outer edges of $G$ and adding just those
inner edges which are in $S$.
(Here the nullity is the number of independent cycles, i.e. the rank
$rk(H_1(G_S, {\mathbb Z}))$ of the first homology group of $G_S$.)
In the formula above
$b_S$ is the unique element of $B_n$ which gives rise to the same partition of the $2n$ boundary points
as $G_S$, in other words $b_S$ is obtained from $G_S$ by contracting each of its inner edges.

The proof that ${\psi}$ is well-defined is analogous to the proof that
the expansion (\ref{chromatic poly2}) of the chromatic polynomial
satisfies the contraction-deletion rule. Specifically, consider the
relations $(1)-(3)$ defining the chromatic algebra (definition
\ref{chromatic definition}). For $(1)$, suppose $e$ is an inner edge
which is not a loop.  The sum (\ref{chromatic expansion}) splits as
the sum over the sets $S'$ which contain the edge $e$ and the sets
$S''$ which do not.  The sets $S'$ are in $1-1$ correspondence with
the sets of edges of $G/e$; the correspondence is given by contracting
$e$. Under this correspondence, the nullity is preserved, and each
term in the sum over the sets $S'$ is equal to the corresponding term
in the sum for the graph $G/e$. Similarly, the sets $S''$ are in $1-1$
correspondence with the sets of edges of $G\smallsetminus e$, and the
corresponding terms in the two sums differ in their sign.  Therefore
${\psi}(G)={\psi}(G/e) - {\psi}(G\smallsetminus e)$, proving the
invariance of ${\psi}$ under the relation $(1)$.

To establish the invariance under $(2)$, let $e$ be a loop in
$G$. Again, (\ref{chromatic expansion}) splits as the sum over the
sets $S'$ which contain $e$ and the sets $S''$ which do not. Both $S'$
and $S''$ are in $1-1$ correspondence with the subsets of inner edges
of $G\smallsetminus e$.  Each term from $S'$ with $e$ has an extra factor
of $Q$ because $n(S')=n(G\smallsetminus e)+1$, while each term from $S''$
has an extra factor of $(-1)$.  Combining these, one gets
${\psi}(G)=(Q-1){\psi}(G\smallsetminus e)$.  The invariance under
$(3)$ is proved analogously: if $G$ contains a $1$-valent vertex, let
$e$ be its adjacent edge. Dividing the subsets into $S'$ and $S''$ as
above, one checks that the corresponding terms cancel in pairs. This
shows that $\psi$ is well-defined.

It follows from the contraction-deletion rule (1) that any graph may be expressed as a linear combination of graphs
without inner edges (elements of $B_n$). Therefore the dimension of ${\mathcal C}_n$ is less than or equal to the cardinality of $|B_n|$.
On the other hand, the
composition $V_n\longrightarrow {\mathcal C}_n\longrightarrow V_n$ is an isomorphism (for $b\in B_n$, ${\psi}({\phi}(b))=\pm b$), proving the opposite inequality
$|B_n|=dim(V_n)\leq dim({\mathcal C}_n)$. Therefore $dim(V_n)=dim({\mathcal C}_n)$ and $B_n$ is a basis of ${\mathcal C}_n$. This completes the proof of lemma \ref{basis lemma}. \qed

\begin{remark}
The seeming difference between the expansions (\ref{chromatic poly2})
and (\ref{chromatic expansion}) is due to the definition of the
chromatic algebra: its trace is the chromatic polynomial of the {\em
dual} graph. The expansion (\ref{chromatic expansion}) is directly
analogous to the expansion for the {\em flow polynomial}, see remark
\ref{Tutte specialization}. Indeed, the quantities $k(S)$ and $n(S)$
in the two expansions correspond to each other under the duality
between $G$ and $\widehat G$.
\end{remark}

With this result at hand, we will now proceed with the proof of
theorem \ref{trivalent lemma}. For each element $\Gamma$ of $B_n$ pick
a trivalent graph $T_{\Gamma}$ such that $H_1(T_{\Gamma})=0$
($T_{\Gamma}$ is a disjoint union of trees) and the contraction of all
inner edges of $T_{\Gamma}$ gives $\Gamma$. Let $TB_n$ denote the set
of such trivalent graphs $\{ T_{\Gamma}, {\Gamma}\in B_n\}$; by
construction $TB_n$ is in a bijective correspondence with $B_n$.  It
is useful for the following argument to note that the combinatorial {\em
F-move}, exchanging the first and third trivalent graphs in figure
\ref{fig:trivalent} and applied at various inner edges of trivalent
graphs, acts transitively on the set of all possible choices of the
graphs $T_{\Gamma}$, for a given ${\Gamma}$.

We claim that any element of ${\mathcal T}_n$ is a linear
combination of elements of $TB_n$. Indeed, given $T\in {\mathcal T}_n$, using the relations in figure \ref{fig:trivalent}
defining ${\mathcal T}_n$, $T$ is seen to be a linear combination of elements $T_i\in {\mathcal T}_n$ where each $T_i$ does not
have cycles: $H_1(T_i)=0$. Now contracting all inner edges of each $T_i$ gives rise to an element ${\Gamma}_i\in B_n$. As noted above,
for each $i$ there is a sequence of $F$-moves taking $T_i$ to the corresponding $T_{{\Gamma}_i}\in TB_n$. Applying the linear relation on
the left in figure \ref{fig:trivalent} each time the $F$-move is needed shows that $T_i=T_{{\Gamma}_i}+$terms with fewer
trivalent vertices. Now an inductive argument, where the induction is on the number of trivalent vertices, proves the claim
that any element of ${\mathcal T}_n$ is a linear
combination of elements of $TB_n$. Using lemma \ref{chromatic expansion}, this shows that
$$dim_{\mathbb C} {\mathcal T}_n\leq |TB_n|=|B_n|=dim_{\mathbb C} {\mathcal C}_n.$$
Combined with the inequality \ref{dim inequality}, this completes the proof of theorem \ref{trivalent lemma}.
\qed

\section{The chromatic algebra and the $SO(3)$ BMW algebra}
\label{sec:equivalence}

In this section we find a homomorphism from the $SO(3)$ BMW algebra to
the chromatic algebra, showing that the defining relations for
$BMW(3)$ reduce to the contraction-deletion rule for planar graphs.
Since the $SO(3)$ Kauffman polynomial arises from the Markov trace of
the $SO(3)$ BMW algebra, this enables us to relate this link invariant
to the chromatic polynomial, and so provide geometric intuition into
the planar description of the algebra.

To give some intuition into why these two algebras are related, it is
useful to recall some results for the integrable field theory
describing the scaling limit of the $Q$-state Potts model near the
critical point. A oft-useful description of an integrable
two-dimensional classical field theory is in terms of the
quasiparticles and their scattering matrices in the corresponding
one-dimensional quantum field theory. For this Potts field theory, the
scattering matrices are invariant under the quantum-group algebra
$U_q(SO(3))$ when the quasiparticles transform in the spin-1
representation \cite{Smir}. This implies that the
scattering matrices can be written in terms of the generators of the
$SO(3)$ BMW algebra.  However, intuitive arguments suggest that the
world lines of the quasiparticles should behave as domain walls
separating regions of like spins \cite{CZ}. These two very different
pictures were reconciled in \cite{FR}. There it was shown how the
generators of the $SO(3)$ BMW algebra subalgebra were those of the
chromatic algebra containing only 4-valent vertices.  This
correspondence was exploited in \cite{FF} to study quantum loop
models. Here we use this motivation to give an elegant
combinatorial-geometric description of $BMW(3)$.

\begin{figure}[ht]
\includegraphics[width=6.5cm]{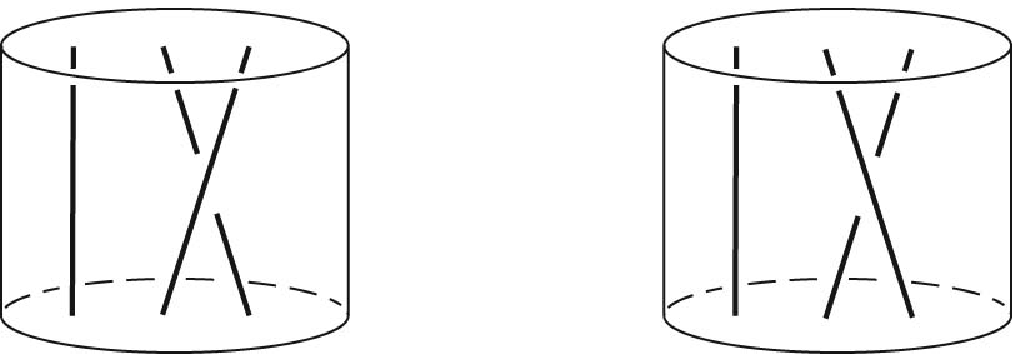}
    \put(-215,28){$B_i =$}
    \put(-100,28){$B^{-1}_i =$}
{\scriptsize
    \put(-157,-7){$i$}
    \put(-142,-7){$i+1$}
    \put(-38,-7){$i$}
    \put(-22,-7){$i+1$}}
    \caption{Braiding generators of $BMW(N)_3$.}
\label{fig:BMW3}
\end{figure}

\begin{figure}[ht]
\includegraphics[height=1.05cm]{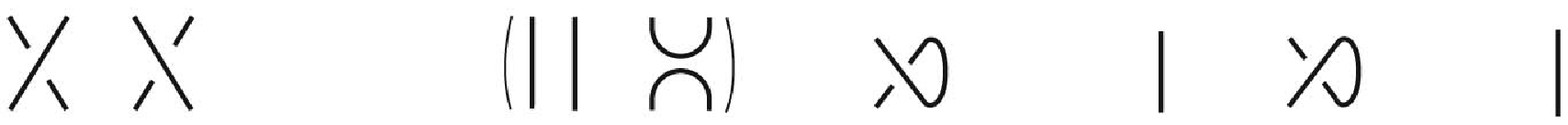}
    \put(-377,12){$-$}
    \put(-338,12){$=\; (q-q^{-1})$}
    \put(-246,12){$-$}
    \put(-207,5){,}
    \put(-148,12){$= q^{1-N}$}
    \put(-99,5){,}
    \put(-46,12){$= q^{N-1}$}
\caption{}
\label{fig:BMWskein}
\end{figure}

We start with a review of the background material on the $SO(N)$
Birman-Murakami-Wenzl algebra; see \cite{BW,M} for
more details.
Here we present the $SO(N)$ BMW algebra in skein form.
A strand can be thought of as corresponding to the fundamental
(dimension $N$) representation of $U_q(SO(N))$.
The braiding
generators (the over/under crossings) are displayed in figure
\ref{fig:BMW3}.
As opposed to $TL$, the skein relations do not allow one to reduce all
braids to non-crossing curves.
Instead, $BMW(N)_n$ is the algebra of framed tangles on $n$
strands in $D^2\times [0,1]$ modulo isotopy and the $SO(N)$ Kauffman
skein relations in figure \ref{fig:BMWskein}.
By a tangle we mean a collection of curves (some of them perhaps
closed) embedded in $D^2\times [0,1]$, with precisely $2n$ endpoints,
$n$ in $D^2\times\{0\}$ and $D^2\times\{1\}$ each, at the prescribed
marked points in the disk. The tangles are framed, i.e. they are given with a
trivialization of their normal bundle. (This is necessary since the last
two relations in figure \ref{fig:BMWskein} are not invariant under the first
Reidemeister move.) As with $TL$, the multiplication is given
by vertical stacking. Like above, $BMW(N)=\cup_n BMW(N)_n$.

We now turn to a presentation of this algebra: the generators of $BMW(N)_n$ include the Temperley-Lieb generators
$1,\, e_1,\ldots,e_{n-1}$ as above, and additionally the braiding
generators $B_i,\, B^{-1}_i$, $i=1,\ldots, n-1$, figure \ref{fig:BMW3}. (This
follows by considering the height Morse function.)  In the algebraic
context, the relations in $BMW_n$ are the Temperley-Lieb relations
(\ref{TL relations}) and in addition
\begin{equation} \label{BMW relations}
B_i e_i=q^{1-N} e_i, \qquad B_i e_{i-1}^{\pm 1} B_i=q^{\pm (N-1)} B_i, \qquad B_i-B_i^{-1}=(q-q^{-1})(1-e_i).
\end{equation}
We will usually work with the geometric
counterparts of these relations: isotopy (expressed as Reidemeister
moves for framed tangles), and the skein relations in figure
\ref{fig:BMWskein}.

Sending the generators of $TL_n$ to the
corresponding generators $e_i$ of
$BMW(N)_n$ defines a map of algebras, and the skein relations require that
deleting a circle has
the effect of multiplying the element of $BMW(N)$ by
\begin{equation}
d_{N}=1+\frac{q^{N-1}-q^{-(N-1)}}{q-q^{-1}}.
\label{dN}
\end{equation}
This gives for example
$d_3=q+1+q^{-1}$ and $d_4=(q+q^{-1})^2$.

The trace, $\hbox{tr}_K\co BMW(N)_n\longrightarrow {\mathbb C}$, is
defined similarly to the TL case.  It is defined on the generators
(framed tangles) by connecting the top and bottom endpoints by
standard arcs in the complement of $D^2\times[0,1]$ in $3$-space,
sweeping from top to bottom, and computing the $SO(N)$ Kauffman
polynomial (given by the skein relations above) of the resulting
link. Below we will discuss this trace in detail.

As with the Temperley-Lieb algebra, it is convenient to distinguish
the skein-theoretic and planar presentations of this algebra,
$BMW(N)^{\rm skein}$ and $BMW(N)^{\rm planar}$.  In the planar version
the additive generators are linear combinations of curves with
crossings, or in other words $4$-valent graphs, in a rectangle.  In
terms of generators, the elements $B, B^{-1}\in BMW^{\rm skein}$ are
replaced with $\crossing\in BMW(N)^{\rm planar}$, defined by \cite{KV}
\begin{equation} \label{basischange}(\crossing)\; =\;
 q(\smoothing)-(\slashoverback)
+q^{-1}(\hsmoothing)\; = \;
 q^{-1}(\smoothing)-(\backoverslash)
+q(\hsmoothing).
\end{equation}
$BMW(N)^{\rm planar}$ is defined as linear combinations of $4$-valent
graphs in a rectangle modulo local relations which are the pull-back
of (\ref{BMW relations}) (or equivalently of the Kauffman skein
relations and the isotopy of tangles) via (\ref{basischange}).  The
equation (\ref{basischange}) yields an isomorphism between the planar
and skein presentations. We will discuss the relations in $BMW^{\rm
planar}$ in more detail below.

The translation of the BMW relations \ref{BMW relations}
(geometrically seen as the skein relations and the isotopy of tangles)
to the planar setting formally follows from the isomorphism
(\ref{basischange}). However, the geometric meaning of the planar
description is not immediately apparent. The results of this paper
provide such a meaning for $BMW(3)$. For the rest of the paper, we will omit the label
and denote this algebra by $BMW$.

One easily checks that there is a map of algebras $TL_n\longrightarrow
{\mathcal C}_n$, where the relation between the Temperley-Lieb and
chromatic parameters is given by $d=Q-1=q+1+q^{-1}$.  This map is
induced by the inclusion $\{$curves$\}\subset\{$graphs$\}$ in a
rectangle $R$. The defining relations for $TL_n$ ($d$-isotopy) hold in
${\mathcal C}_n$ due to the relation (2) in definition \ref{chromatic
definition} of the chromatic algebra. The following statement shows
that this map extends to $BMW_n$.

\begin{them} \label{equivalence theorem} \sl
The formulas
\begin{equation} \label{BMW to chromatic}
(\slashoverback) \,  \mapsto  \,
 q(\smoothing)-(\crossing)
+q^{-1}(\hsmoothing), \;\;\;
(\backoverslash)\,   \mapsto  \,
 q^{-1}(\smoothing)-(\crossing)
+q(\hsmoothing)
\end{equation}
define a homomorphism of algebras  $i\co BMW^{\rm skein}_n\longrightarrow {\mathcal C}_n$ over ${\mathbb C}[q]$.
\end{them}

\begin{remark}
It follows from theorem \ref{relations theorem} below that the homomorphism $i$ preserves the traces
of these algebras, see corollary \ref{Kauffman chromatic Corollary}.
\end{remark}

\begin{remark} \label{planar remark}
Abusing the notation, we will keep the same symbol for the map $i\co
BMW^{\rm planar}_n\longrightarrow {\mathcal C}_n$.  This homomorphism
is induced by the inclusion $\{4$-valent graphs$\}\subset \{$all
graphs$\}$. Recall that the relations in $BMW^{\rm planar}_n$ are the
pullback of (\ref{BMW relations}) under the identification
(\ref{basischange}), and the content of the theorem above is that
these relations are a consequence of the contraction-deletion rule.
\end{remark}

\begin{figure}[ht]
\includegraphics[height=1.7cm]{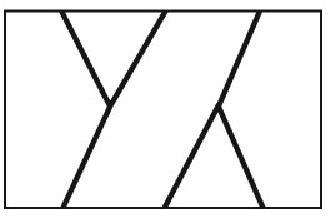}
\caption{}
\label{fig:image}
\end{figure}

\begin{remark} \label{BMW remark}
The homomorphism $i$ is not surjective (one may check that for example, the graph in figure \ref{fig:image} is not equivalent in ${\mathcal C}_3$ to a linear combination of $4$-valent graphs), however it seems reasonable to conjecture that it is injective, in other words that $BMW_n$
is a subalgebra of the chromatic algebra ${\mathcal C}_n$ generated by $4$-valent graphs. Note that the chromatic algebra
${\mathcal C}_n$ has quite elegant presentations in terms of the contraction-deletion rule (definition \ref{chromatic definition}), and in terms of
trivalent graphs (definition \ref{trivalent definition}). In contrast, the corresponding presentation in terms of $4$-valent graphs
is rather involved: the relations are the pullback of (\ref{BMW relations}) under the identification (\ref{basischange}). See also
a related discussion in section \ref{BMW planar chromatic}.
\end{remark}

\begin{figure}[ht]
\includegraphics[width=14cm]{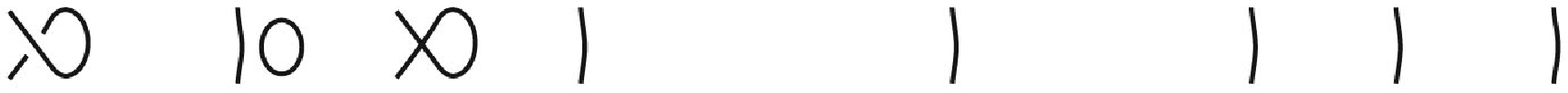}
{\scriptsize
    \put(-372,10){$=\, q^{-1}$}
    \put(-314,10){$-$}
    \put(-272,10){$+\, q$}
    \put(-243,10){$=\, q^{-1}(q+1+q^{-1})$}
    \put(-150,10){$-\, (q+1+q^{-1})$}
    \put(-71,10){$+\, q$}
    \put(-37,10){$=\,  q^{-2}$}}
     \caption{}
\label{fig:proof1}
\end{figure}

{\em Proof of theorem \ref{equivalence theorem}.} One needs to show that the map $i$ is well-defined. The first defining relation of the BMW algebra (figure \ref{fig:BMWskein})
$$(\slashoverback)-(\backoverslash)\, =\, (q-q^{-1})[(\smoothing)-
(\hsmoothing)]$$ follows directly from the equations (\ref{BMW to
chromatic}). One needs to check that the last two relations in figure
\ref{fig:BMWskein}, as well as the regular isotopy of tangle diagrams
(the second and third Reidemeister moves), hold in the chromatic
algebra.  The proof of the first of these is shown in figure
\ref{fig:proof1}.

\begin{figure}[ht]
\includegraphics[height=.9cm]{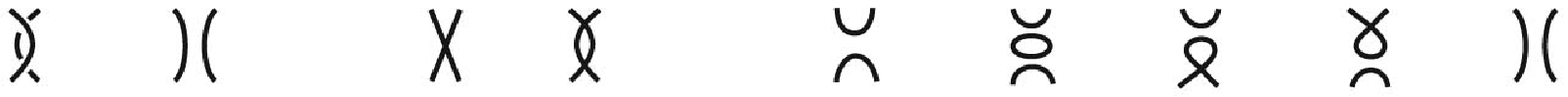}
{\scriptsize
    \put(-385,11){$=$}
    \put(-349,11){$-\,(q+q^{-1})$}
    \put(-279,11){$+$}
    \put(-248,11){$+\,(q^2+q^{-2})$}
    \put(-168,11){$+$}
    \put(-127,11){$-\,q$}
    \put(-86,11){$-\,q^{-1}$}
    \put(-35,11){$=$}}
     \caption{}
\label{fig:Reidem2}
\end{figure}

The last relation in figure \ref{fig:BMWskein} is established analogously.
To prove the second Reidemeister move, start with the diagram on the left in figure \ref{fig:Reidem2} and resolve the crossings according to the formulas (\ref{BMW to chromatic}). Applying the contraction-deletion rule (1) in definition \ref{chromatic definition} of the chromatic algebra to the edges connecting the double points in the third term, eliminating
the trivial circle in the last term according to the rule (2), and canceling the resulting terms, one gets the diagram on the right.

\begin{figure}[ht]
\includegraphics[height=.9cm]{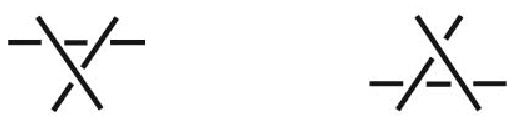}
     \caption{}
\label{fig:Reidem3}
\end{figure}

\begin{figure}[ht]
\includegraphics[height=.9cm]{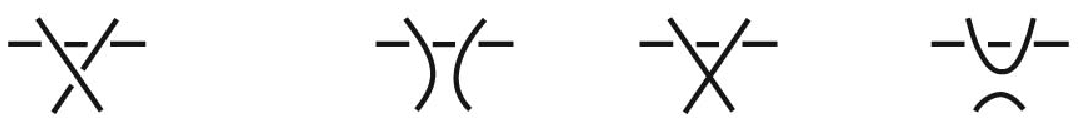}
{\scriptsize
    \put(-200,10){$=\;\;\, q^{-1}$}
    \put(-117,10){$-$}
    \put(-57,10){$+\; q$}}
     \caption{}
\label{fig:Reidem31}
\end{figure}

The remaining relation is the third Reidemeister move: one has to show that the images of the two diagrams in
figure \ref{fig:Reidem3} are equal in the chromatic algebra. Note that these two diagrams differ by a 180 degree rotation, so
whenever a planar diagram, invariant under such a rotation, appears in the expansion of one of them, it also
appears (with the same coefficient)
in the expansion of the other one. Expanding the lower crossing  of the diagram on the left according to (\ref{BMW to chromatic}),
one gets the expression in figure \ref{fig:Reidem31}.

It follows from the remark above, and from the proof of the second Reidemeister move
that the first and the third terms in figure \ref{fig:Reidem31} cancel with the corresponding terms in the expansion
of the second diagram in figure \ref{fig:Reidem3}. Therefore it remains to show that the second term on the right in figure \ref{fig:Reidem31} equals
its 180 degree rotation in the chromatic algebra.

\begin{figure}[ht]
\includegraphics[height=.8cm]{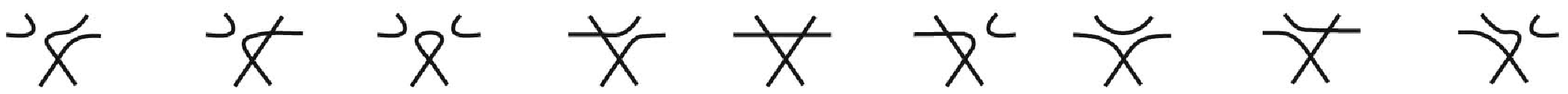}
{\scriptsize
    \put(-420,10){$q^{-2}$}
    \put(-376,10){$-q^{-1}$}
    \put(-320,10){$+$}
    \put(-280,10){$-q^{-1}$}
    \put(-228,10){$+$}
    \put(-186,10){$-q$}
    \put(-140,10){$+$}
    \put(-97,10){$-q$}
    \put(-50,10){$+q^2$}}
     \caption{}
\label{fig:Reidem32}
\end{figure}

\begin{figure}[ht]
\includegraphics[height=1cm]{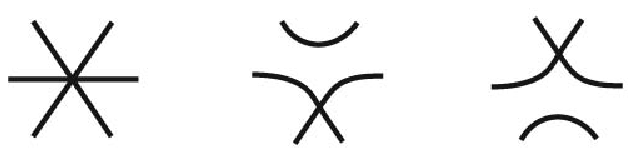}
{\scriptsize
    \put(-182,12){$(q-1+q^{-1})$}
    \put(-90,12){$+$}
    \put(-43,12){$+$}}
     \caption{}
\label{fig:Reidem33}
\end{figure}

The expansion of this term according to (\ref{BMW to chromatic}) is shown in figure \ref{fig:Reidem32}. The 4th and 8th terms
are invariant under the 180 rotation. Omitting these two terms, and using the relations (1), (2) in the chromatic
algebra to expand the terms with more than one double point, one gets the expression in figure \ref{fig:Reidem33}.

This expression is again invariant under a 180 degree rotation, and this concludes the proof
for the third Reidemeister move and the proof of theorem \ref{equivalence theorem}.
\qed

\section{Relations between $TL$, $BMW$, and the chromatic algebra}
\label{sec:relations}

In this section we investigate the relations between the $SO(3)$ Birman-Murakami-Wenzl, chromatic, and
Temperley-Lieb algebras. The main result is stated in theorem \ref{relations theorem}. This relationship
is useful in a variety of contexts. For example, it allows one to use the Jones-Wenzl projectors
in $TL^d$, at special values of $d$, to analyze the structure of the trace radical of the BMW and chromatic
algebras. This was used in \cite{FK2} to find linear relations obeyed by the chromatic polynomial of planar
graphs, evaluated at Beraha numbers. This also provides a relationship between certain string-net models and loop
models, cf \cite{FFNWW}, \cite{Fnew}. Applying traces to these algebras, we express the $SO(3)$ Kauffman polynomial
in terms of the chromatic polynomial, see corollary \ref{Kauffman chromatic Corollary}.

The $SO(3)$ BMW algebra may be described in several ways.
First we explain how each strand of the BMW algebra may be viewed as a ``fusion'' of two
Temperley-Lieb strands, in other words defining a homomorphism of the $SO(3)$ BMW algebra in degree $n$,
$BMW_n$ to $TL_{2n}$. This description is particularly
natural when studying quantum-group algebras, where each strand in $TL$
and $BMW$ correspond respectively to a spin-1/2 and spin-1
representation of $U_q(sl_2)$.

In fact, we define a homomorphism from the chromatic algebra ${\mathcal C}_n\longrightarrow TL_{2n}$
which yields a map from the BMW algebra by pre-composing it with the homomorphism $i\co BMW_n\longrightarrow {\mathcal C}_n$,
constructed in the previous section. Recall that the chromatic algebra ${\mathcal C}_n$ is defined
as the quotient of the free algebra ${\mathcal F}_n$ by the ideal generated by the relations $(1)-(3)$ in
definition \ref{chromatic definition}.

\begin{figure}[ht]
\includegraphics[width=11.5cm]{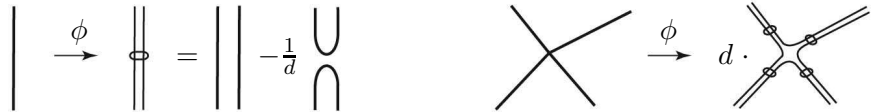}
    \put(-305,29){${\phi}$}
    \put(-265,19.5){$=$}
    \put(-234,19.5){$-\frac{1}{d}$}
    \put(-82,29){${\phi}$}
    \put(-60,19.5){$d\;\cdot$}
\caption{Definition of the homomorphism ${\phi}\co {\mathcal C}^Q_n\longrightarrow TL^d_{2n}$, where $Q=d^2$.}
\label{fig:phi}
\end{figure}

\begin{defi}    \label{def:phi}
Define a homomorphism ${\phi}\co {\mathcal F}_n\longrightarrow TL_{2n}$ on the additive generators
(graphs in a rectangle) of the free graph algebra ${\mathcal F}_n$
by replacing each edge with the linear combination ${\phi}(\;|\;)=\smoothing-\frac{1}{d}\,\hsmoothing\,$, and resolving
each vertex as shown in figure \ref{fig:phi}.
The factor in the definition of ${\phi}$ corresponding to a
$r$-valent vertex is $d^{(r-2)/2}$, so for example it equals $d$ for the
$4$-valent vertex in figure \ref{fig:phi}. The overall factor
for a graph $G$ is the product of the factors $d^{(r(V)-2)/2}$ over
all vertices $V$ of $G$.
\end{defi}

\begin{figure}[ht]
\includegraphics[height=1.6cm]{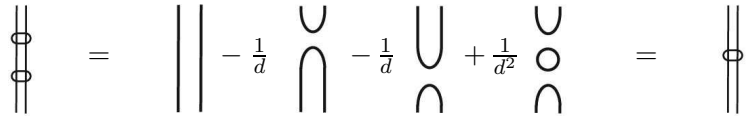}
    \put(-250,21){$=$}
    \put(-200,21){$-\;\frac{1}{d}$}
    \put(-151,21){$-\,\frac{1}{d}$}
    \put(-108,21){$+\,\frac{1}{d^2}$}
    \put(-43,21){$=$}
\caption{$P_2\circ P_2=P_2$}
\label{fig:JW}
\end{figure}

Therefore ${\phi}(G)$ is a sum of $2^{E(G)}$ terms, where $E(G)$ is the number of edges of $G$.
Note that ${\phi}$ replaces each edge with the
second Jones-Wenzl projector $P_2$, well-known in the study of the Temperley-Lieb algebra \cite{J}.
They are idempotents: $P_2\circ P_2=P_2$, and this identity (used in the definition of the homomorphism
$\phi$) may be easily checked directly, figure \ref{fig:JW}.

Various authors have considered versions of the map ${\phi}$ in the
knot-theoretic and TQFT contexts, see \cite{Y,Ja,KL, FFNWW,Walker}. In \cite{FR}
this was used to give a map of the SO(3) BMW algebra to
the Temperley-Lieb algebra.

\begin{lem} \sl ${\phi}$ induces a well-defined homomorphism of algebras ${\mathcal C}^Q_n\longrightarrow TL^d_{2n}$, where
$Q=d^2$.
\end{lem}

\begin{figure}[ht]
\includegraphics[height=1.8cm]{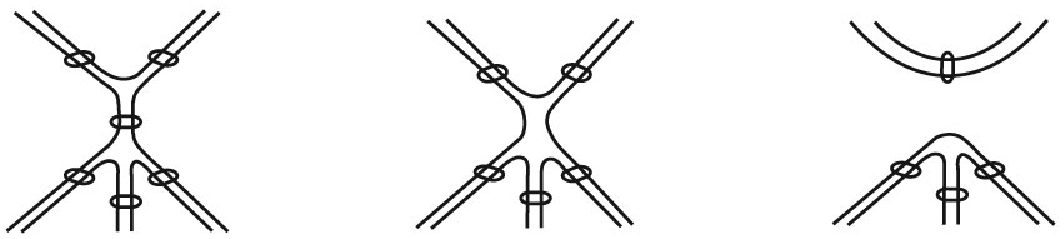}
    \put(-168,22){$=$}
    \put(-76,22){$-$}
\caption{}
\label{fig:phi3}
\end{figure}

{\em Proof.} One needs to check that ${\phi}$ is well-defined with respect to the relations $(1)-(3)$
in the chromatic algebra (definition \ref{chromatic definition}).
To establish (1), one applies ${\phi}$ to both sides and expands
the projector at the edge $e$, as shown in figure
\ref{fig:phi}. The resulting relation holds due to the choices
of the powers of $d$ corresponding to the valencies of the
vertices, figure \ref{fig:phi3}.

Similarly, one uses the definition to check the relations (2) and (3).
\qed

To state the main result of this section, we also need to define a homomorphism
${\phi}'\co BMW_n^{\rm skein}\longrightarrow TL_{2n}^{\rm skein}$ (the latter algebra is defined by (\ref{Jones relation}))
where the relation between the
parameters $q$ in the $BMW$ algebra and $A$ in $TL$ is given by $q=A^4$. This map is the
``$2$-coloring'': it is
defined on tangle generators of $BMW_n^{\rm skein}$ by replacing each strand with
the Jones-Wenzl projector $P_2$. One may check that ${\phi}'$ is well-defined directly from definitions
(also see \cite[p.35]{KL}), and this also follows from the commutativity of the diagram below.

\begin{them} \label{relations theorem}
\sl The following diagram commutes.
\begin{equation}    \label{commutative diagram}
\xymatrix{ {\mathbb C} \ar[rrrr]^{=}  \ar[ddd]_{=} & & & & {\mathbb C} \ar[ddd]^{=} \\
    &   BMW^{\rm skein}_n \ar[ul]_{tr_{K}}  \ar[d]_{\cong} \ar[rr]^{{\phi}'}  & &
            TL^{\rm skein}_{2n} \ar[ur]^{tr_{\langle\rangle}} \ar[d]^{\cong} &  \\
    &   BMW^{\rm planar}_n \ar[dl]_{tr}  \ar[r]^{\;\;\;\;\;\;\;\;i} &
    {\mathcal C}_n \ar[d]_{tr_{\chi}} \ar[r]^{\!\!\!\!\!\!\!\!\!\!{\phi}} & TL^{\rm planar}_{2n} \ar[dr]^{tr_d} & \\
    {\mathbb C} \ar[rr]^{=}  & & {\mathbb C} \ar[rr]^{=} & & {\mathbb C}
}
\end{equation}
\end{them}

For convenience of the reader, we recall the notations in the diagram above. The parameter $q$ in the $BMW$
algebras, $Q$ in the chromatic algebra ${\mathcal C}_n$, $A$ in $TL^{\rm skein}_{2n}$ and $d$ in $TL^{\rm planar}_{2n}$ are
related by:
$$ q=A^4, \;\; d=-A^2-A^{-2}, \;\; Q=q+2+q^{-1}=d^2.$$
The traces of various algebras are defined on their respective generators as follows:

$tr_K\co BMW_n^{\rm skein}\longrightarrow {\mathbb C}$ is given by the $SO(3)$
Kauffman polynomial, figure \ref{fig:BMWskein}.

$tr\co BMW_n^{\rm planar}\longrightarrow {\mathbb C}$ is the pull-back of $tr_K$ via the
isomorphism (\ref{basischange}): $BMW_n^{\rm planar}\cong BMW_n^{\rm skein}$.

$tr_{\chi}\co {\mathcal C}_n\longrightarrow {\mathbb C}$ is $Q^{-1}$ times the chromatic
polynomial of the dual graph, see (\ref{chromatic trace}).

$tr_{\langle\rangle}\co TL_{2n}^{\rm skein}\longrightarrow {\mathbb C}$ is given by the Kauffman bracket defined
in (\ref{Jones relation}).

$tr_d\co TL_{2n}^{\rm planar}\longrightarrow {\mathbb C}$, discussed
in the section \ref{sec:TL}, is computed as $d^{\# loops}.$

\begin{figure}[ht]
\includegraphics[height=1.8cm]{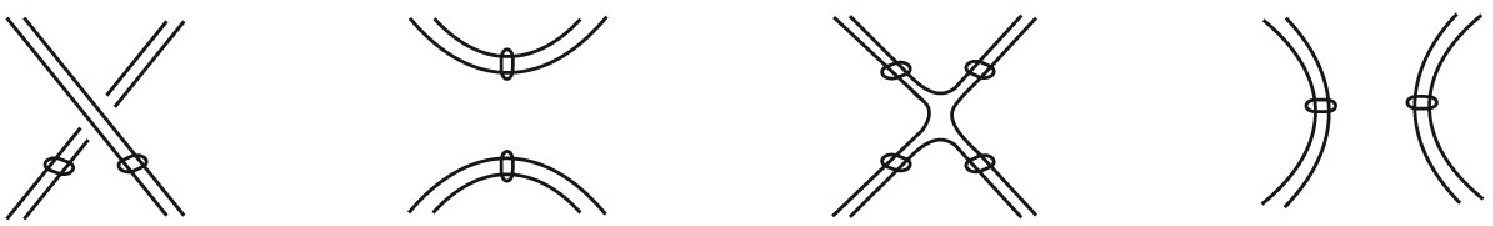}
    \put(-290,24){$=\;\;\; q\,\cdot$}
    \put(-190,24){$-\;\;\; d\,\cdot$}
    \put(-90,24){$+\;\; q^{-1}\,\cdot$}
\caption{}
\label{fig:phi4}
\end{figure}

{\em Proof of theorem \ref{relations theorem}.} We begin the proof by showing that the interior diagram
of algebra homomorphisms (without the traces) commutes. This amounts to showing that the two maps in the
diagram from $BMW^{\rm skein}_n$ to $TL^{\rm planar}_{2n}$ send the braiding generator to the
same element. The braiding generator $(\backoverslash)$ is mapped by ${\phi}'$ to the element of $TL$ shown on the
left in figure \ref{fig:phi4}. (The crossings are resolved using (\ref{Jones relation}) to get an element of $TL^{\rm planar}_{2n}$.)
The other map, $BMW^{\rm skein}\cong BMW^{\rm planar}_n \longrightarrow {\mathcal C}_n \longrightarrow TL^{\rm planar}_{2n}$,
sends $(\backoverslash)$ to the linear combination shown on the right in figure \ref{fig:phi4}. (The middle term
acquires the coefficient $d$ since the vertex is $4$-valent, see definition \ref{def:phi}.)

This identity in the Temperley-Lieb algebra is established in \cite[p.35]{KL}. We now consider the traces in
the diagram (\ref{commutative diagram}).

\begin{lem} \label{commutative} \sl Let $G$ be a planar graph. Then
  for $Q=d^2$,
\begin{equation} \label{chromatic dual}
tr_\chi(G)=Q^{-1}\, {\chi}_Q(\widehat G)\, =\, tr_d ({\phi}(G)).\end{equation}
Therefore, the following diagram
commutes:
\begin{equation} \label{chromaticTL}
\xymatrix{ {\mathcal C}_n^{\scriptsize Q}  \ar[d]^{tr_{\chi}} \ar[r]^{\phi} & TL_{2n}^{d}  \ar[d]^{tr_{d}}\\
{\mathbb C}\ar[r]^=  &  {\mathbb C} }\end{equation}
\end{lem}

{\em Proof.} The proof (involving the expansion (\ref{chromatic poly2}) of the chromatic polynomial)
for trivalent graphs $G$
is given in lemma 2.5 in \cite{FK2}. Using the contraction-deletion rule (1) in definition \ref{chromatic definition},
any graph may be represented as a linear combination of trivalent graphs. The statement then follows from the
fact that map $\phi$ and the traces in the diagram above are well-defined. \qed

Observe that the identity in figure \ref{fig:phi4} shows that the map ${\phi}'\co BMW^{\rm skein}_n\longrightarrow TL^{\rm skein}_{2n}$
preserves the traces: the evaluation of the $SO(3)$ Kauffman polynomial of the closure of a framed tangle equals the evaluation
of the Kauffman bracket of the $2$-coloring of that link. Recall that the isomorphism $BMW_n^{\rm planar}\cong BMW_n^{\rm skein}$
preserves the trace by definition. Examining the diagram (\ref{commutative diagram}), one observes that the remaining map
$i\co BMW^{\rm planar}_n \longrightarrow {\mathcal C}_n$ preserves the traces as well. This concludes the proof of
theorem \ref{relations theorem}. \qed

We state the last observation in the proof above as a corollary. Given
a crossing $(\backoverslash)$, one calls $(\hsmoothing)$ its
$0$-resolution, and $(\smoothing)$ its $1$-resolution. Given a planar
diagram $D$ of a framed link $L$, use (\ref{basischange}) to express
it as a linear combination of planar $4$-valent graphs $G_i$. For
each $i$, let $p_i$ denote the number of $0$-resolutions and $n_i$ the
number of $1$-resolutions that are used in (\ref{basischange}) to get
the graph $G_i$ from the diagram $D$. Let $v(G_i)$ denote the number
of vertices of $G_i$.

\begin{cor} \label{Kauffman chromatic Corollary} \cite{Ja} Given a framed link $L$, using the notations above
the $SO(3)$ Kauffman polynomial $K_L$ may be expressed as
$K_L(q)\; =\; Q^{-1} \sum (-1)^{v(G_i)}\, q^{p_i-n_i}\, {\chi}_{\widehat{G_i}}(Q)$, where the summation is taken over
all $4$-valent planar graphs $G_i$ which are the result of applying the formula (\ref{basischange}) at each crossing of
a planar diagram
of $L$. Here $Q=(q^{1/2}+q^{-1/2})^2$, and ${\chi}_{\widehat{G_i}}$ denotes the chromatic polynomial of the graph dual to
$G_i$.
\end{cor}

\section{Properties of the chromatic inner product} \label{sec:inner product}

We use the homomorphism ${\phi}$ to the Temperley-Lieb algebra, defined in
\ref{def:phi}, to
study the trace pairing structure on the chromatic algebra ${\mathcal C}_n$
(introduced in \ref{chromatic trace def}.)

\begin{lem} \label{lem:injective} \sl
The homomorphism ${\phi}\co {\mathcal C}_n^{\scriptsize {d^2}}
 \longrightarrow TL_{2n}^{d}$ is injective.
\end{lem}

\begin{figure}[ht]
\includegraphics[height=2cm]{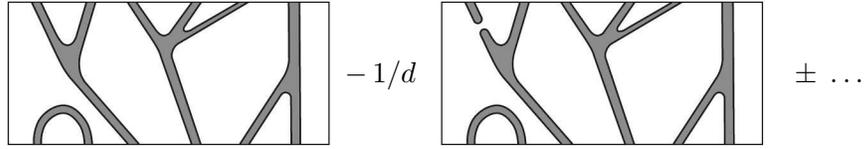}
\put(-160,25){$-\, 1/d$}
\put(10,25){$\pm \; \ldots$}
\caption{${\phi}$ applied to the graph in figure \ref{fig:basis}.}
\label{fig:injective}
\end{figure}

{\em Proof.}
It is easy to see that for any graph $G$ in the rectangle $R$, the terms in
the expansion ${\phi}(G)$ are
in $1-1$ correspondence with the subsets $S\subset E(G)$ of the set of edges
of $G$.
More precisely, given such $S$, consider the graph $G_S$ whose vertex set
consists of all
vertices of $G$ and the edge set is $S$. Then the terms of ${\phi}(G)$
correspond
(with the coefficient depending on $|S|$) to the boundary
of the regular neighborhood of the graphs $G_S$. (For an edge $e$, the two
terms in the definition of $\phi$ in figure \ref{fig:phi} correspond to the
two possibilities: $e\in S$, $e\not\in S$.)
For more details see the proof of lemma 2.5 in \cite{FK2}.

In lemma \ref{basis lemma} we showed that (the isotopy classes of) graphs
without inner edges in the rectangle $R$ form a basis of ${\mathcal C}_n$.
Suppose
a linear combination of such graphs, $\sum {\alpha}_i G_i$, is in the kernel
of
${\phi}$.
Consider ${\phi}(G_i)$ for some $i$, see figure \ref{fig:injective} which
corresponds to the graph in figure \ref{fig:basis}. Consider the ``leading
term''
of ${\phi}(G_i)$, where each edge of $G_i$ is replaced by its two parallel
copies
(the first term in the example in figure \ref{fig:injective}.) In terms of
the
correspondence discussed above, this term is given by $S=E(G_i)$.
We claim that
this term does not arise in the expansion of any other ${\phi}(G_j),$ $j\neq
i$.

It is clear that it does not arise as a leading term for any other graph,
since each graph
is recovered as the spine of its leading term. Since the graphs $G_j$ do not
have inner edges, it
also does not arise as a non-leading term
of any ${\phi}(G_j)$, since any such term involves a ``shaded turn-back'' --
corresponding to
a $G_S$ which has an isolated vertex on the boundary of $R$, as in the
second term in figure \ref{fig:injective}.
Such ``shaded turn-backs'' cannot be part of a leading term for
${\phi}(G_i)$ since the graphs
with an isolated vertex in the interior of $R$ are trivial in the chromatic
algebra (relation $(3)$ in
definition \ref{chromatic definition}.) This concludes the proof of lemma
\ref{lem:injective}.
\qed

\begin{cor} \label{cor:inner product} \sl Given $Q\geq 0$, the pairing ${\mathcal C}^Q\otimes
{\mathcal C}^Q\longrightarrow {\mathbb C}$,
introduced in \ref{chromatic trace def}, is a positive definite Hermitian
product provided that $Q\geq 4$.
For $Q$ equal to Beraha numbers $B_n=2+2\cos(2{\pi}/(n+1))$, the pairing is
positive semidefinite.
\end{cor}

The corollary follows from the corresponding statement for the
Temperley-Lieb algebra \cite{J} at $d=\sqrt Q$ and lemmas
\ref{commutative} and \ref{lem:injective}.

We note the intriguing fact that the values of $Q$ for which the Hermitian
product on
${\mathcal C}^Q_n$ is positive-definite coincides with the values of $Q$
($Q\geq 4$)
for which the chromatic polynomial of planar graphs is conjectured to be
positive
\cite{BL}. (Of course, $Q=4$ corresponds to the $4$-color theorem.) At
Beraha numbers, the
trace pairing descends to a positive-definite product on the quotient of
${\mathcal C}$ by
the trace radical, and it gives rise to the unitary structure of the
(doubled) $SO(3)$ TQFT,
see section \ref{TQFT section}.

Among the elementary consequences of the corollary above, consider the
Cauchy-Schwarz inequality: given two graphs $G, H$ in a disk with an equal
number of endpoints on the boundary of the disk, for $Q\geq 4$, $(G,H)^2\leq
(G,G) \cdot (H,H)$.
The inner product is given by gluing two disks along their boundary and
computing the chromatic polynomial
at $Q\geq 4$ of the dual graph of the resulting graph in the $2$-sphere.
This may be loosely
formulated as saying that the chromatic polynomial at these values of $Q$ of
such graphs is maximized by graphs with
a reflection symmetry, a statement that does not seem to follow immediately
from the combinatorial definition
of the chromatic polynomial.

\section{Concluding remarks and questions} \label{sec:questions} We mention a number of questions
motivated by our results.

\subsection{} \label{chromatic question}
The chromatic algebra ${\mathcal C}$ defined in this paper provides a
natural framework for studying algebraic and combinatorial properties
of the chromatic polynomial. We mentioned some of the elementary
consequences in sections \ref{sec:chromatic}, \ref{sec:trivalent}. In
\cite{FK2} we use this algebra to give an algebraic proof of Tutte's
golden identity \cite{T} and to establish chromatic polynomial
relations evaluated at Beraha numbers, whose existence was conjectured
by Tutte.  The chromatic algebra should be useful for a variety of
other problems as well, for example it seems likely that it should
give new insight into Birkhoff-Lewis equations \cite{BL} and the
associated Tutte's invariants \cite{T1}. (See also \cite{CJ} for more
details on this subject.  The use of the Temperley-Lieb algebra in
\cite{CJ} appears to be different from our approach; it would be
interesting to find a connection with our results.)

\subsection{}
The chromatic algebra and its relations with the TL and BMW algebras should also be useful in analyzing analytic properties of the
chromatic polynomial. Specifically, Tutte established the estimate $|{\chi}^{}_T({\phi}+1)|\leq
{\phi}^{5-k}$, where $T$ is a planar triangulation and
$k$ is the number of its vertices. The value ${\phi}+1$ is one of Beraha numbers, $B_5$, where
$B_n=2+2\cos(2{\pi}/n)$. Considering the map
${\phi}\co {\mathcal C}_n\longrightarrow TL_{2n}$ (defined in section
\ref{sec:relations}), the Beraha numbers correspond to the special values of $d$, $d_n=2+2\cos({\pi}/n)$,
where the trace radical of the Temperley-Lieb algebra is non-trivial, and is generated by
the Jones-Wenzl projector \cite{J}. This algebraic structure may prove useful in determining
whether there is an analogue of Tutte's estimate for other Beraha numbers. (This question is
interesting in connection with the observation that the real roots of the chromatic polynomial
of large planar triangulations seem to accumulate near the points $\{ B_n\}$.)

In sections \ref{sec:equivalence}, \ref{sec:relations} we established a relationship between
the $SO(3)$ BMW and chromatic algebras. The former is useful for studying the quantum $SO(3)$
invariants of $3$-manifolds (via their surgery presentation), while the latter is directly
related to the chromatic polynomial.
It is known \cite{LW}, \cite{Wong} that the $SO(3)$
quantum invariants of $3$-manifolds are dense in the complex plane.
This motivates the question of whether there may be a related
density result for the
values of the chromatic polynomial of planar graphs at Beraha numbers.

\subsection{} \label{BMW planar chromatic} The $SO(N)$ Kauffman polynomial may be viewed as an invariant of
planar $4$-valent graphs via the ``change of basis'' formula (\ref{basischange}). (Similarly the generators of the
$SO(N)$ BMW algebra may be taken to be $4$-valent planar graphs, rather than tangles.) The relations among $4$-valent
graphs are then the pullback of the Kauffman relations in figure \ref{fig:BMWskein} and of the last two Reidemeister moves
under (\ref{basischange}).

In this paper we show that in the case $N=3$ all of these relations follow from the chromatic relations in definition
\ref{chromatic definition}. We also discuss the special case
$N=4$ in \cite[section 4]{FK2}, relating $BMW(4)$ with $TL\times TL$.
It is an interesting question whether there is a nice combinatorial/geometric
interpretation of the $SO(N)$ invariant of graphs for other values of
$N$ involving unlabeled graphs. (A different approach is taken in
\cite{Ku} in the rank $2$ case where the graph edges have different
labels.)

\subsection{} In this paper we define and investigate some of the properties of the map
from the $SO(3)$ BMW algebra to the chromatic algebra, $i\co
BMW_n\longrightarrow {\mathcal C}_n$, see theorem \ref{equivalence theorem}.  It would be interesting to
establish a precise relationship between these algebras.  As indicated
in remark \ref{BMW remark}, we believe this homomorphism is injective
but not surjective.  (From the TQFT perspective, an interesting
question is whether these two algebras are Morita equivalent.)  It
follows from lemma \ref{commutative} that at special values of
$Q=2+2\cos(2{\pi}j/n)$, the pull-back of the Jones-Wenzl projector
from the Temperley-Lieb algebra via the homomorphism $\phi$ is in the
trace radical of the chromatic algebra. A question important from the
perspective of finding linear relations among the values of the
chromatic polynomial of planar graphs (see \cite{FK2}) is whether this
pull-back of the Jones-Wenzl projector generates the entire trace
radical of the chromatic algebra.


\smallskip









\begin{thebibliography}{10}

\bibitem{Baxbook} R.J.~Baxter,
{\em Exactly Solved Models in Statistical Mechanics} (Academic,
London, 1982)

\bibitem{B} B. Bollobas, {\em Modern graph theory}, Springer, 1998.

\bibitem{BL} G.D. Birkhoff and D.C. Lewis, {\em Chromatic polynomials}, Trans. Amer. Math. Soc. 60 (1946), 355-451.

\bibitem{BW} J. Birman and H. Wenzl, {\em Braids, link polynomials and a new
algebra}, Trans. Amer. Math. Soc. 313 (1989), 249--273.

\bibitem{CJ} S. Cautis and D. Jackson, {\em The matrix of chromatic joins and the Temperley-Lieb algebra},
J. Combin. Theory Ser. B 89 (2003), 109--155.

\bibitem{CZ}
L.~Chim and A.~B.~Zamolodchikov,
{\em Integrable field theory of q state Potts model with $0 < q < 4$},
Int.~J.~Mod.~Phys A 7 (1992) 5317.

\bibitem{Fnew}
P.~Fendley,
{\em Topological order from quantum loops and nets}, Ann. Physics 323 (2008), 3113-3136. [arXiv:0804.0625]

\bibitem{FF}
P.~Fendley and E.~Fradkin,
{\em Realizing non-Abelian statistics in time-reversal-invariant
  systems},
Phys.~Rev.~B 72 (2005) 024412 [arXiv:cond-mat/0502071]

\bibitem{FJ} P.~Fendley and J.L.~Jacobsen, {\em Critical points in
  coupled Potts models and critical phases in coupled loop models},
J. Phys. A: Math. Theor. 41 (2008) 215001
[arXiv:0803.2618]

\bibitem{FK2}
P.~Fendley and V.~Krushkal,
{\em Tutte chromatic identities from the Temperley-Lieb algebra},
Geom. Topol. 13 (2009), 709-741 [arXiv:0711.0016]

\bibitem{FR}
P.~Fendley and N.~Read,
{\em Exact S-matrices for supersymmetric sigma models and the Potts
  model},
J.~Phys. A 35 (2003) 10675
[arXiv:hep-th/0207176]


\bibitem{FFNWW}L. Fidkowski, M. Freedman, C. Nayak, K. Walker, Z. Wang,
{\em From String Nets to Nonabelions},
cond-mat/0610583.

\bibitem{FK}
C.M.~Fortuin and P.W.~Kasteleyn,
{\em On The Random Cluster Model.~1.~Introduction And Relation To Other Models},
Physica {\bf 57}, 536 (1972).


\bibitem{F} M. Freedman, {\em A magnetic model with a possible Chern-Simons
phase},  Comm. Math. Phys.  234  (2003), 129--183 [arXiv:quant-ph/0110060]

\bibitem{FNWW} M. Freedman, C. Nayak, K. Walker, Z. Wang, {\em  On Picture (2+1)-TQFTs},
arXiv:0806.1926.

\bibitem{Ja} F.~Jaeger, {\em On some graph invariants related to the Kauffman polynomial},
Progress in knot theory and related topics. Paris: Hermann. Trav. Cours. 56, 69-82 (1997).

\bibitem{J} V.F.R.~Jones,
{\em Index for subfactors},  Invent. Math. 72 (1983), 1-25.

\bibitem{J1} V.F.R.~Jones,
{\em Planar algebras}, arXiv:math/9909027.

\bibitem{KL} L.H. Kauffman and S.L. Lins,
{\em Temperley-Lieb recoupling theory and invariants of $3$-manifolds},
Princeton University Press, Princeton, NJ, 1994.

\bibitem{KV} L.H. Kauffman and P. Vogel, {\em Link polynomials and a graphical
calculus}, J. Knot Theory Ramifications  1  (1992), 59--104.

\bibitem{KooSaleur} W.M. Koo and H. Saleur, 
{\em Fused Potts models},
Internat. J. Modern Phys. A 8 (1993), 5165--5233. 


\bibitem{Ku} G. Kuperberg, {\em Spiders for rank $2$ Lie algebras},
 Comm. Math. Phys.  180  (1996), 109--151.

\bibitem{Larsen} M. Larsen and Z. Wang, {\em Density of the SO(3) TQFT representation of mapping class groups},
 Comm. Math. Phys.  260  (2005), 641--658.

\bibitem{LW}
M.~A.~Levin and X.~G.~Wen, {\em String-net condensation: A physical mechanism for topological phases}
Phys.~Rev.~B 71 (2005) 045110
[arXiv: cond-mat/0404617]

\bibitem{MW} P. P. Martin and D. Woodcock, {\em The partition algebras and a new
deformation of the Schur algebras}, J. Algebra 203 (1998), 91-124.

\bibitem{M} J. Murakami, {\em The Kauffman polynomial of links and representation theory.}
Osaka J. Math. 24 (1987), 745--758.

\bibitem{RT} N. Reshetikhin, V. Turaev, {\em Invariants of $3$-manifolds via link polynomials and quantum
groups},  Invent. Math.  103  (1991), 547-597.

\bibitem{Smir}
F.A.~Smirnov,
{\em Exact S matrices for phi(1,2) perturbated minimal models of
  conformal field theory}
Int.~J.~Mod.~Phys.~A 6 (1991) 1407.

\bibitem{TL}
H.~N.~V.~Temperley and E.~H.~Lieb,
{\em Relations between the 'percolation' and 'colouring' problem and other
graph-theoretical problems associated with regular planar lattices: some
exact results for the 'percolation' problem},
  Proc.\ Roy.\ Soc.\ Lond.\  A {\bf 322}, 251 (1971).

\bibitem{Turaev} V. Turaev, Quantum invariants of knots and 3-manifolds,
Walter de Gruyter \& Co., Berlin, 1994.

\bibitem{TV} V. Turaev and O. Viro, {\em State sum invariants of $3-$manifolds and
quantum $6j-$symbols}, Topology 31 (1992), 865-902.


\bibitem{T0} W.T. Tutte, {\em On chromatic polynomials and the golden ratio},
J. Combinatorial Theory 9 (1970), 289--296.

\bibitem{T} W.T.~Tutte, {\em More about chromatic polynomials and the golden ratio}, Combinatorial Structures
and their Applications, 439-453 (Proc. Calgary Internat. Conf., Calgary, Alta., 1969)

\bibitem{T1} W.T. Tutte, {\em On the Birkhoff-Lewis equations}, Discrete Math., 92 (1991), 417-425.


\bibitem{Walker} K.~Walker, talk given at IPAM conference
  ``Topological Quantum Computing'', available at
{\tt https://www.ipam.ucla.edu/schedule.aspx?pc=tqc2007}

\bibitem{W} K. Walker,
  {\em On Witten's $3$-manifold invariants}, available at
  http://canyon23.net/math/

\bibitem{Witten}   E.~Witten,
  {\em Quantum field theory and the Jones polynomial},
  Commun.\ Math.\ Phys.\  {121} (1989) 351.

\bibitem{Wong} H. Wong, {\em $SO(3)$ quantum invariants are dense}, preprint.

\bibitem{Y} S. Yamada, {\em An operator on regular isotopy invariants of link diagrams},  Topology  28
(1989), 369--377.







\end{thebibliography}
\end{document}